\magnification=\magstep1%metka
\input amstex
\documentstyle{amsppt}
\TagsOnRight
\hsize=5.1in                                                  
\vsize=7.8in

\topmatter

\title Von Neumann categories \\ 
and extended $L^2$ cohomology
\endtitle
%\rightheadtext{}
%\leftheadtext{}
\author  Michael Farber \endauthor
\address
School of Mathematical Sciences,
Tel-Aviv University,
Ramat-Aviv 69978, Israel
\endaddress
\email farber\@math.tau.ac.il
\endemail
\thanks{The research was supported by a grant from the
US - Israel Binational Science Foundation}
\endthanks
\abstract{In this paper we suggest a new general formalism for studying the
$L^2$ invariants of polyhedra and manifolds. First, we examine generality 
in which one may apply the
construction of {\it the extended abelian category}, which was suggested 
in \cite{Fa, Fa1}, using the ideas of P.Freyd \cite{F}. 
This leads to the notions of {\it a finite von Neumann category} 
and of {\it a trace} on such category. Given a finite von Neumann category, 
we study {\it the
extended $L^2$ homology and cohomology} theories with values in the abelian
extension. Any trace on the initial category produces numerical invariants -
the von Neumann dimension and the Novikov - Shubin numbers. Thus, we obtain 
{\it the local versions of the Novikov - Shubin invariants}, localized 
at different 
traces. In the {\it "abelian"} case this localization can be made more
geometric: we show that any torsion object determines a
{\it "divisor"}  -- a closed subspace of the space of the parameters.
The divisors of torsion objects together with the information produced by 
the local Novikov - Shubin invariants may be used to study {\it multiplicities of 
intersections} of algebraic and analytic varieties (we discuss 
here only simple 
examples demonstrating this possibility). We compute explicitly the divisors 
and the von Neumann dimensions of the extended $L^2$ cohomology in the real 
analytic situation. We also give general formulae for the extended $L^2$ 
cohomology of a mapping torus. Finally, we show how one can define a De Rham
version of the extended cohomology and prove a De Rham type theorem.}
\endabstract
\endtopmatter
%-------------------------------------------------------------------------
%----------------------------------------------------------------------

\define\C{{\Bbb C}}
\define\R{{\Bbb R}}  
       
\define\Z{{\Bbb Z}}

\define\OO{{\Cal O}}  
\define\T{{\Cal T}}  
\define\X{{\Cal X}}   
\define\Y{{\Cal Y}}  
\define\M{{\Cal M}}

\define\Hom{\operatorname{Hom}}  
\define\om{{\frak {hom}}}

\define\Tr{\operatorname{Tr}}

\define\Tot{\operatorname{Tot}}
   
\define\Ext{\operatorname{Ext}}

\define\End{\operatorname{End}}

\define\A{{\Cal A}}  
\redefine\B{{\Cal B}}  
\define\e{{\frak e}} 
\redefine\c{{\frak c}}

\define\tr{\operatorname{tr}} 
\define\im{\operatorname{im}} 
   
%\define\Ker{\operatorname{Ker}}   

\define\id{\operatorname{id}}    
\define\Id{\operatorname{Id}}
\define\ns{\frak {ns}}     
\define\ob{\operatorname{ob}}   
\define\cl{\operatorname{cl}} 
\define\supp{\operatorname{supp}}  
\define\rank{\operatorname{rank}}  
%\define\ln{\operatorname{ln}}
\redefine\H{\Cal H}
\redefine\D{\Cal D}
\define\E{\Cal E}
\define\U{{\Cal U}}  
\def\<{\langle}
\def\>{\rangle}

\define\pd#1#2{\dfrac{\partial#1}{\partial#2}}
%\define\fa{\Cal F(A)}
   
\define\cc{\Cal C}
  
\define\eca{\Cal E({\Cal C}_{\Cal A})}

\define\ca{\Cal C_{{\Cal A}}}

\define\cabp{\Cal C_{{\Cal A}\otimes{\Cal B'}}}   
\define\tca{\T(\ca)}
\redefine\L{\Cal L} 
\define\ts{\tilde\otimes_\Lambda}

\define\caf{\ca^{fin}}
\define\fa{\Cal F_{\A}}

\documentstyle{amsppt}   

\nopagenumbers

\vskip 2cm

The classical approach to the $L^2$ cohomology theory, as developed in 
\cite{A}, \cite{CG}, \cite{D}, \cite{G}, consists in viewing 
it as a functor which assigns Hilbert modules over a von 
Neumann algebra to polyhedra and manifolds.
This functor is also called {\it the reduced $L^2$ 
cohomology}, since in order to preserve the Hilbert module structure
on the cohomology space one defines this cohomology as the quotient of 
the space
of cocycles modulo the {\it closure} of the space of coboundaries.

This approach has certain well-known problems caused by the fact that 
Hilbert modules over a von Neumann algebra form only an {\it additive 
category, and not an abelian category}. The obtained cohomology
theory has, in particular, difficulties in
dealing with exact sequences. On the other hand, the reduced $L^2$ cohomology
theory does not feel the phenomena {\it "zero in the continuous spectrum"},
discovered by S.Novikov and M.Shubin \cite{NS}, \cite{NS1}, which carries
interesting topological information. This phenomenon
is also related to the problem of absence of exactness; it
was recently studied in \cite{GS}, \cite{GS1} and in \cite{LL}.

In  \cite{Fa}, \cite{Fa1} there were constructed new homology and cohomology 
theories with values in {\it the extended abelian category}, containing 
(as the full subcategory of projectives) the usual category of finitely 
generated Hilbertian modules over a finite von Neumann algebra. 
These theories are called {\it extended $L^2$ homology
and cohomology}. The {\it reduced} $L^2$ (co)homology 
appears as {\it the projective part} of the extended (co)homology. 
The extended homology and cohomology contain
also {\it the torsion parts}; it was shown in \cite{Fa, Fa1} that the 
torsion parts of the extended $L^2$ cohomology determine {\it the 
Novikov-Shubin invariants}. Moreover, it turned out (cf. \cite{Fa, Fa1}) 
that the torsion parts of the extended cohomology produce also some 
new numerical invariants,
which are independent of the Novikov-Shubin number. With the aid of these
new invariants the Morse type inequalities of S.Novikov and 
M.Shubin \cite{NS, NS1} were strengthened in \cite{Fa, Fa1}, so that the
spectra near zero of the Laplacians also give quantitative information
about the critical points.

The construction of the extended abelian category, which was suggested 
in \cite{Fa, Fa1}, used the general categorical study of P.Freyd \cite{Fr}.
It was described in \cite{Fa, Fa1} in the simplest and smallest possible
version, extending the
additive category of finitely generated modules over a von Neumann algebra
supplied with a fixed finite trace. In this small version the extension
construction is equivalent to the purely algebraic alternative construction
which was suggested later by W. L\"uck \cite{L}.

One of the main goals of the present paper is to
investigate this extension construction in full generality. We try to
find requirements on the category of Hilbert representations of a given
$\ast$-algebra under which the extension construction works. This leads to
the notion of {\it a Hilbert category}, which is assumed to be closed
under taking kernels and
adjoint morphisms. We show that any Hilbert category can be canonically 
embedded in an abelian category. We establish a similar result for the 
Hilbertian representations as well, cf. \S 5. However, it is clear that
it is impossible to apply this construction in the case of Frechet spaces
or Banach spaces.

The further development of the theory
can be made under additional assumption that the initial category is
{\it a finite von Neumann category}. This guarantees that 
the torsion subcategory is also 
an abelian category. We show that any von Neumann algebra 
determines a von Neumann category; we discuss many interesting 
examples of von Neumann categories.

I should mention that (as I learned when
the present paper was completely finished), 
P. Ghez, R. Lima, and J.E. Roberts in 1985 studied the notion of $W^\ast$-category
(cf. \cite{GLR}),
which is equivalent to the notion of von Neumann category of this paper.
In \cite{GLR} they proved that many important results of the theory of von
Neumann algebras (including, for example, the modular theory)
can be generalized to the setting of von Neumann categories.
Their motivation in \cite{GLR}
was completely different from ours, and this fact clearly explains why the
intersection between the present paper and \cite{GLR} is so little.

Next notion which appears to be extremely important for the formalism
developed here is the notion of {\it a trace on a von Neumann category}.
By the definition, a trace on a von Neumann category, is a compatible family
of traces (in the usual sense of von Neumann algebras) on the rings of 
endomorphisms of all the objects. For example, the category of finite
dimensional vector spaces has only one trace up to normalization.
Traces allow us {\it to localize} the notions of {\it von Neumann dimension} 
and the {\it Novikov - Shubin invariants}. 
The idea is to fix a category and to vary
the trace on it. Given a trace on a finite von Neumann category, 
it determines the von Neumann dimension function (measuring sizes of the 
projective objects) and also the {\it Novikov - Shubin number} (measuring 
sizes of the torsion objects). We establish here the main properties of these
function, which are mainly known, although in a different contexts.

As a special interesting example we consider here (in \S 4)
the {\it abelian case}, 
which is given by the von Neumann  category generated by fields of 
Hilbert spaces over a fixed compact space $Z$. This category allows to 
study {\it families}. Any chain complex formed by a
sequence of vector bundles over $Z$ and bundle maps between them determines
a complex in this von Neumann category by considering the spaces of $L^2$
sections. We show that one may gain some
interesting geometric information by studying the extended cohomology of the
obtained complexes of $L^2$-sections. The projective part of the extended
cohomology (i.e. the reduced $L^2$ cohomology) can be easily explicitly 
computed. The torsion part of the extended cohomology contains much more
information. In particular, it determines a {\it divisor} - a closed subset 
of $Z$. We prove
that in the real analytic situation the divisor coincides with the subset 
of point $\xi\in Z$ where the fiberwise cohomology is not generic.

Simple examples computed in \S 4 show that one may hope to be able to 
study the {\it multiplicities of intersections} of algebraic and analytic 
varieties using the technique of extended cohomology.

The functors of extended homology and cohomology on the category of finite 
polyhedra are defined in section
\S 6. As a coefficient systems for these theories we use arbitrary modules
over the group ring in the extended abelian category (allowing torsion
objects). This construction
generalizes the well established practice of considering the regular 
representation only. We point out some basic properties of this construction,
and then as an application, we compute explicitly the extended cohomology
of mapping tori.

In the last section \S 7 we prove a De Rham type theorem for the extended
cohomology. There are two candidates for the De Rham complex -
one which is built on infinitely smooth forms and the other, which uses
Sobolev spaces. We first show that the smooth De Rham complex is homotopy 
equivalent to the combinatorially defined \v Cech complex. However,
one cannot use the smooth De Rham complex to define extended cohomology,
since it is not clear how to construct an abelian category containing
this complex. On the contrary, the Sobolev - De Rham complex 
(we show that it is homotopy equivalent to the smooth De Rham complex) 
belongs to the category of Hilbertian representations, and so the abelian 
extension construction of \S 5 applies to it.

I should add that M. Shubin was the first who initiated discussions of De Rham type
theorems for the extended cohomology. 
His recent preprint \cite{S} contains a De Rham type theorem different from
the theorems of \S 7 of the present paper. I am very thankful to M. Shubin who
sent to me a very preliminary version of \cite{S}.

I would like to mention also the preprint \cite{CCMP}, which discusses
problems of similar nature, concerning
relations between the spectral properties of smooth Laplacian and its
combinatorial analogues corresponding to different triangulations.

I am very thankful to A. Carey, M. Gromov, P. Milman, H. Moscovichi and M. Shubin 
for useful discussions and help.

\heading{\bf \S 1. Abelian extensions of Hilbert categories}\endheading 

In this section we will describe a construction of an 
abelian category extending in a canonical way a given category of 
$\ast$-representations of an algebra with involution. 
It is based on the general categorical study of P.Freyd \cite{F}.
This construction appears to be of fundamental importance for the homology 
theory of 
topological algebras and their Hilbert representations. A special case of 
this construction was used in \cite{Fa, Fa1} to introduce the notion of
extended $L^2$ cohomology and to explain the homological nature of the 
Novikov-Shubin invariants.

\subheading{1.1} Let $\A$ be an algebra over $\C$ having an involution 
which will be denoted 
by the star $\ast$. Recall, that this means that $\ast:\A\to\A$ is an 
involutive skew-linear map such that $(ab)^\ast=b^\ast a^\ast$ for 
$a,b\in \A$. {\it A Hilbert representation of $\A$ (or a Hilbert module)} 
is a Hilbert space $\H$ supplied
with a left action of $\A$ on $\H$ by bounded linear maps such that for any
$a\in \A$ holds 
$$\<ax,y\>=\<x,a^\ast y\>\tag1-1$$ 
for all $x,y\in \H$. {\it A morphism} between Hilbert representations 
$\phi: \H_1\to \H_2$ is a bounded linear map commuting with the action of 
the algebra $\A$. 

The category of Hilbert representations is an additive category.
It is easy to see that the kernel of any morphism $\phi:\H_1\to \H_2$ 
between Hilbert representations in the sense of category theory
coincides with $\phi^{-1}(0)$. However one observes that the categorical 
image of any morphism $\phi:\H_1\to \H_2$ coincides with the {\it closure}
of the set-theoretic image $\cl(\phi(\H_1))$. 

The category of Hilbert representations 
has the following important property: if $\H_1\subset\H_2$ 
is a {\it closed} $\A$-invariant submodule then the orthogonal complement of 
$\H_1$ in $\H_2$ (denoted by $\H_1^\perp$) is also
a closed $\A$ submodule and the orthogonal projection $\pi:\H_2\to \H_1$ is
an $\A$-homomorphism. Thus, every closed submodule has a complement.
Equivalently, any surjective morphism  
$\phi:\H_1\to \H_2$ between Hilbert representations splits.

Next we observe that {\it the category of Hilbert representations is not an 
abelian category.} Indeed, any morphism $\phi:\H_1\to \H_2$ between Hilbert 
representations which is injective and has dense image is both monomorphic 
and epimorphic; it is not an isomorphism unless $\phi$ is onto. 

\subheading{1.2} Our aim in this section is to construct, following 
ideas of P.Freyd \cite{F}, a larger {\it abelian} category, containing 
the category of Hilbert representations as a full subcategory.

In applications, we will always have some additional structures on the 
algebra $\A$
(for example, topology, trace, etc.) and we will restrict ourselves to 
special additive subcategories of the category of Hilbert 
$\A$-representations. 
Our aim is to describe the construction of the extended abelian category 
in the most general unified way, in order to include all applications. 

With this goal in mind,
we will assume that an additive subcategory $\cc_\A$ of the category of 
Hilbert representations over $\A$ is specified. We will suppose 
that this additive subcategory $\cc_\A$ has the following properties:
\roster 
\item"{(i)}" {\it The kernel of any morphism $\phi:\H_1\to \H_2$ in
$\cc_\A$  and the natural inclusion $\ker \phi\to \H_1$ belong to $\cc_\A$.
\item"{(ii)}" For any morphism $\phi:\H_1\to \H_2$ of $\ca$ the adjoint
operator $\phi^\ast:\H_2\to \H_1$ is also a morphism of $\ca$.}
\endroster
Observe, that from (i) and (ii) follow:
\roster
\item"{(iii)}" {\it The closure of the image $\cl(\im\phi)$ of any morphism 
$\phi:\H_1\to \H_2$ in $\cc_\A$ and also the natural projection 
$\H_2\to \H_2/\cl(\im\phi)$ belong to $\cc_\A$.
\item"{(iv)}" Suppose that $\H_1\subset\H_2$ is a closed $\A$-invariant 
submodule. 
If $\H_1$,
$\H_2$ and the inclusion $\H_1\to \H_2$ belong to $\cc_\A$ then the orthogonal
complement $\H_1^\perp$ and the inclusion $\H_1^\perp\to \H_2$ belong to
$\cc_\A$.}
\endroster

\proclaim{Definition} An additive category $\ca$ of $\ast$-representations
of a $\ast$-algebra $\A$ satisfying the conditions
(i) and (ii) above will be called {\it Hilbert category}.\endproclaim

We will postpone discussion of examples until section 2.6.

In applications to topology, we will consider the situation 
when $\A$ is the group 
algebra of a discrete group $\C[\pi]$. 

\subheading{1.3. Abelian extension of a Hilbert category} Returning to 
the general situation, given a Hilbert category 
$\ca$, we are going to define a bigger 
category $\E(\cc_\A)$, containing $\ca$ as a full subcategory. 

{\it An object} of the category $\E(\cc_\A)$ is defined as a
morphism $(\alpha:A^\prime\to A)$ in the category $\cc_\A$. 
Given a pair of objects $\X=(\alpha: A^\prime \to A)$ and 
$\Y=(\beta:B^\prime\to B)$ of $\E(\cc_\A)$,
a {\it morphism} $\X\to\Y$ in the category $\E(\cc_\A)$ is an  
equivalence class of morphisms $f:A\to B$ of category $\cc_\A$
such that $f\circ\alpha=\beta\circ g$
for some morphism $g:A^\prime \to B^\prime$ in $\cc_\A$. 
Two morphisms $f:A\to B$ and $f^\prime:A\to B$ of $\cc_\A$ represent 
{\it identical morphisms $\X\to\Y$ of $\E(\cc_\A)$} iff 
$f-f^\prime = \beta\circ F$ for some morphism $F:A\to B^\prime$ of category
$\cc_\A$. This defines an equivalence relation.
The morphism $\X\to\Y$, represented by $f:A\to B$,
is denoted by
$$[f]:(\alpha:A^\prime\to A)\ \to\ (\beta:B^\prime\to B)\quad
\text{or by}\quad
[f]:\X\to\Y.$$
The {\it composition} of morphisms is defined as the composition of the 
corresponding morphisms $f$ in the category $\cc_\A$.

The category $\E(\cc_\A)$ will be called {\it the abelian extension of 
the category} $\cc_\A$ or {\it the extended category}, for short. We will see
in Proposition 1.7 below that it is indeed abelian.

\subheading{1.4. Excision} The following construction (which we will call 
{\it excision}) describes the 
typical changes which sometimes may be performed on objects of the extended 
category $\E(\cc_\A)$ without changing their isomorphism class.

Suppose that $\X=(\alpha: A^\prime \to A)$ is an object of $\E(\cc_\A)$ 
and $P\subset A^\prime$ is a closed $\A$ submodule which belongs to
$\cc_\A$ and such that its image $Q=\alpha(P)\subset A$ is also closed. 
Using the assumptions (i)-(ii) concerning the category $\cc_\A$ we obtain
the following object $\Y=(\beta:P^\perp\to Q^\perp)$ where $\beta$ is the
composition of $\alpha|_{P^\perp}$ and the orthogonal projection onto 
$Q^\perp$. We leave it as an easy exercise to show that $\Y$ is isomorphic 
to the original object $\X$ inside $\E(\cc_\A)$.

Using the excision, {\it we may represent any isomorphism type in $\eca$
by an injective morphism} $(\alpha:A^\prime\to A)$ (since we may always 
make an excision with respect to the kernel).

\subheading{1.5} Here is another simple observation. {\it An object 
$(\alpha:A^\prime\to A)$
of $\eca$ represents the zero object in $\eca$ (i.e. is isomorphic to 
$(0\to 0)$) if and only if $\alpha$ is surjective.} In fact, if $\alpha$ is 
surjective
then we may make an excision with respect to whole $A^\prime$ which produces
the zero object. On the other hand, if $(\alpha:A^\prime\to A)$ is isomorphic
to the zero object, then from the definitions we obtain that there exists
$F:A\to A^\prime$ with $\alpha\circ F = 0$; thus $\alpha$ is surjective.

In order to show that $\E(\cc_\A)$ {\it is an abelian category} we will need
the following statement.

\proclaim{1.6. Proposition} Let
$[f]:(\alpha:A^\prime\to A)\ \to\ (\beta:B^\prime\to B)$          
be a morphism in $\E(\cc_\A)$.
Then its kernel is represented by
$$[k]:\ (\gamma:P^\prime \to P)\ \to\ (\alpha:A^\prime\to A),$$
where
$$
\botaligned
{\CD
P@>k>>A\\
@V{f^\prime}VV  @VVfV\\ 
B^\prime@>\beta>>B&
\endCD}
\qquad\text{{and}}\qquad
{\CD
P^\prime@>{k^\prime}>>A^\prime\\
@V{f^{\prime\prime}}VV  @VV{f\circ \alpha}V\\ 
B^\prime@>\beta>>B&
\endCD}
\endbotaligned\tag1-2
$$
are the pullbacks of the diagrams

$$
\botaligned
{\CD
  @. A\\
@.   @VVfV\\ 
B^\prime@>\beta>>B
\endCD}
{\qquad\text{and}\qquad}
{\CD
  @. A^\prime\\
@.   @VV{f\circ\alpha}V\\ 
B^\prime@>\beta>>B
\endCD}
\endbotaligned\tag1-3
$$
correspondingly, and $\gamma:P^\prime\to P$ is the canonical map, induced 
by the obvious map of the right diagram (1-2) into the left one.

The cokernel of the above morphism $[f]$
%$$[f]:(\alpha:A^\prime\to A)\ \to\ (\beta:B^\prime\to B),$$  
is represented by
$$[\id_B]:\ (\beta:B^\prime\to B)\ \to\ ((\beta,-f):B^\prime\oplus A\to B).$$
\endproclaim

\demo{Proof} Note that pullbacks are defined as follows: $P$
is the kernel of the morphism $f\oplus{-\beta} : A\oplus B^\prime\to B$,
and $P^\prime$ is defined similarly; they exist by our assumption (i) above. 

To prove the statement concerning the kernels, note first that clearly
$[f]\circ [k]=0$. Suppose that $(\gamma:C^\prime\to C)$ is an abject of 
$\E(\cc_\A)$ and that we are given a morphism 
$$[g]:(\gamma:C^\prime\to C)\to (\alpha:A^\prime\to A)$$
such that $[f]\circ [g]=0$. We want to show that $[g]$ could be factorized
uniquely through $[k]$. Since $f\circ g$ represents a zero morphism in
$\eca$, there exists a morphism $F:C\to B^\prime$ such that 
$\beta\circ F= f\circ g$. The pair of morphisms $(g,F)$ determines a morphism
$h:C\to P$ with $k\circ h = g$. From the definition of morphisms in $\eca$
we know that there is a morphism 
$g^\prime: C^\prime\to A^\prime$ with $\alpha\circ g^\prime = g\circ\gamma$;
thus, the pair of morphisms $(g^\prime,F\circ\gamma)$ determines a morphism
$h^\prime:C^\prime\to P^\prime$ showing that we have a morphism is $\eca$
$$[h]: (\gamma:C^\prime\to C)\to (k:P^\prime\to P)$$
which gives the desired factorization $[g]=[k]\circ[h]$.

To show uniqueness of the above factorization, we may assume in the notations
of the previous paragraph that $[k]\circ [g]=0$ and then prove that then 
$[g]=0$. In fact, $[k]\circ [g]=0$ means (according to our definitions) 
that there exists a morphism $G:C\to A^\prime$ with $k\circ G=\alpha\circ G$. 
The pair of morphisms $G:C\to A^\prime$ and $C@>g>>P\to B^\prime$ determine 
a morphism $G^\prime:C\to P^\prime$ with $\gamma\circ G^\prime = g$; 
thus $[g]=0$.

The second statement concerning the cokernels can be checked similarly.
\enddemo

\proclaim{1.7. Proposition} $\eca$ is an abelian category.
\endproclaim
\demo{Proof} We already know that any morphism in $\eca$ has a kernel and a
cokernel. We have to prove that any morphism 
$[f]:(\alpha:A^\prime\to A)\ \to\ (\beta:B^\prime\to B)$  in $\eca$ 
which is both monomorphic and epimorphic is an isomorphism in $\eca$. We may 
assume without loss of generality, that both $\alpha$ and $\beta$ are
injective. We will use the notations introduced in Proposition 1.6. The 
morphism $\gamma:P^\prime\to P$ is an isomorphism and we have a exact
sequences 
$$0\to P@>{(k,f^\prime)}>>A\oplus B^\prime@>{(f,-\beta)}>>B\to 0$$
$$0\to P^\prime @>{(k^\prime,f^{\prime\prime})}>>A^\prime\oplus B^\prime
@>{(f^\prime,-\id)}>>B^\prime\to 0$$      
Thus we may represent the morphism $[f]$ as the composite of two excisions
$$(\alpha:A^\prime\to A)@>{\simeq}>>(
\bmatrix \alpha & 0\\0 & 1\endbmatrix:
A^\prime\oplus B^\prime\to A\oplus B^\prime)
@>{(f,-\beta)}>{\simeq}>(\beta:B^\prime\to B)
$$
Here the excision on the right is performed with respect to the closed 
submodule $P^\prime\subset A^\prime\oplus B^\prime$; it is mapped by
$\bmatrix \alpha & 0\\ 0 & 1\endbmatrix$ onto the closed submodule
$P\subset A\oplus B^\prime$. 

This completes the proof. $\square$\enddemo

\subheading{1.8. Embedding of $\ca$ into $\eca$} Given an object $A\in\ob(\ca)$
one defines the following object $(0\to A)\in\ob(\eca)$ of the extended
category. Since any morphism $f:A\to B$ determines a morphism 
$[f]:(0\to A)\to (0\to B)$ in the extended category, we obtain an
embedding $\ca\to \eca$. In fact, it is easy to see that this embedding is
full.

We want to characterize the objects of the extended category which are
isomorphic in $\eca$ to objects coming from $\ca$ in intrinsic terms:

\proclaim{1.9. Proposition} An object $\X\in\ob(\eca)$ is projective if and
only if it is isomorphic in $\eca$ to an object of the form $(0\to A)$,
where $A\in\ob(\ca)$.\endproclaim
\demo{Proof} Let $\X=(\beta:B'\to B)$; by 1.4 we may assume that $\beta$ is
injective. Using Proposition 1.6 we see that the sequence
$$
0\to (0\to B')@>{[\beta]}>>(0\to B)@>{[\id]}>>\X\to 0
$$
is exact and so it splits if $\X$ is supposed to be projective. Any splitting
is given by a homomorphism $r:B\to B$ such that $r\circ \beta = 0$ and
$\id - r = \beta\circ F$ for some $F:B\to B'$. It now follows that the
image of $\beta$ coincides with the kernel of $r$ and so the image of $\beta$
is closed. Thus, we may make excision with respect to the whole $B'$ and hence
by 1.4 we obtain that $\X$ is isomorphic in $\eca$ to $(0\to A)$, where
$A = B/\beta(B')$.

Conversely, let $\X=(0\to A)$ and suppose that 
$[f]:(\alpha:A'\to A)\to (\beta:b'\to B)$ is an
epimorphism in $\eca$. We want to show that arbitrary morphism 
$[g]: (0\to A)\to (\beta:B'\to B)$ can be lifted in $(\alpha:C'\to C)$.
By Proposition 1.6, the morphism 
$(\beta,-f):B'\oplus C\to B$ is an epimorphism in $\ca$ and so there exists
a morphism $G:A\to B'\oplus C$ with $(-\beta,f)\circ G = g$. Writing
$G = g_1\oplus g_2$, where $g_1:A\to B'$ and $g_2:A\to C$, we obtain
$[g] = [f]\circ [g_2]$ is $\eca$.
$\square$.\enddemo

Using embedding $\ca\to\eca$ and Proposition 1.9, any chain complex in $\ca$ 
can be viewed as {\it a projective chain complex in category} 
$\eca$. Since $\eca$ is abelian, we may define the homology of the chain 
complex as a graded object of
$\eca$. The corresponding homology will be called {\it the extended homology
of the chain complex}. The following statement contains the explicit
computation of the extended homology.

\proclaim{1.10. Proposition} Given a chain complex $C$ in a Hilbert category 
$\ca$
$$\dots \to C_{i+1}@>{\partial}>>C_i
@>{\partial}>>C_{i-1}@>{\partial}>>\dots,\tag1-4$$
the homology of this complex, viewed as an object of the extended category 
$\eca$, equals
$$H_i(C) = (\partial: C_{i+1}\to Z_i) \quad \in \quad \ob(\eca) ,\tag1-5$$
where $Z_i$ is the space of the cycles, 
$Z_i = \ker[\partial: C_i\to C_{i-1}]$.
\endproclaim
\demo{Proof} It is an easy exercise based on applying Proposition 1.6.
\enddemo

The term {\it extended homology} intends to emphasize distinction between
the present approach and the well established practice of assigning the
following (sometimes it is called {\it reduced}) homology to a chain 
complex of Hilbert spaces
$$Z_i/\cl(\im [\partial C_{i+1}\to C_i]).$$ 
The extended homology is {\it larger}
since there is natural epimorphism in $\eca$
$$H_i(C) = (\partial: C_{i+1}\to Z_i) \to 
Z_i/\cl(\im [\partial C_{i+1}\to C_i]).$$

\subheading{1.11. Functor of projective part} Here we will introduce the 
notion which allows to express the relation between the reduced
and extended homology in a precise way.

Let $\X=(\alpha:A'\to A)$ be an object of the extended category $\eca$.
Its {\it projective part} is defined as the following object of $\ca$
$$P(\X) = A/\cl(\im(\alpha)).\tag1-6$$
Clearly, any morphism $f:\X\to\Y$ in $\eca$ induces a morphism in $\ca$
$P(f):P(\X)\to P(\Y)$ between the projective parts and so {\it we have 
well defined functor} 
$$P: \eca \to \ca.\tag1-7$$

For a short exact sequence in $\eca$
$$0\to \X' \to \X \to \X'' \to 0,$$
the corresponding sequence of projective parts
$$0\to P(X') \to P(\X) \to P(\X'') \to 0\tag1-8$$
may be {\it not exact} in the middle term (although it is exact in the 
other places). Thus the functor of projective part {\it is not half-exact}.

Using the functor of projective part we may express the rule which explains
how the extended homology determines the reduced homology. Namely,
{\it the reduced homology
coincides with the projective part of the extended homology}.

Non-exactness of the functor $P$ of projective part explains
well-known difficulties with exact sequences of the reduced homology. 
The extended homology theory is free of these problems since it assumes
its values in an honest abelian category.

\subheading{1.12. Remark} Suppose that a Hilbert category $\ca$ is a subcategory
of another Hilbert category $\ca'$. Then we have an obvious functor
$F:\eca\to \E({\Cal C}'_{\A})$ between the extended categories. In general,
it may happen that $F([f])=0$ for a nonzero morphism $[f]$ of $\eca$.

However, {\it if a Hilbert category $\ca$ is a full subcategory of another
Hilbert category $\ca'$ then the above functor $F$ represents the extended
abelian category $\eca$ as a full subcategory of the extended category
$\E({\Cal C}'_{\A})$.}

\heading {\bf \S 2. Von Neumann categories}\endheading

Here we will introduce a notion of von Neumann category which naturally
generalizes the notion of von Neumann algebra.
The notion of von Neumann category is equivalent (but it is not identical)
to the notion of
$W^\ast$-category which was first introduced in \cite{GLR}.
It seems possible to 
develop the theory of von Neumann categories to a large extend in a way 
similar to the classical theory of von Neumann algebras. We will however 
mention here only basic properties, which are used in this paper.

Our interest in von Neumann categories is based on the fact, that
(as we will see later in section \S 3) the extended abelian category $\eca$
has many important additional properties in case when the original category
$\ca$ is a von Neumann category.

\subheading{2.1. Von Neumann categories of Hilbert representations} 
We will say that a category of Hilbert representations $\ca$ of a 
$\ast$-algebra $\A$ is {\it a von Neumann category} if it satisfies besides
conditions (i) and (ii) of section 1, the following condition:

\roster  
\item"{(v)}" {\it for any pair of representations $\H_1, \H_2\in\ob(\ca)$, 
the corresponding 
set of morphisms $\Hom_{\ca}(\H_1,\H_2)$ is a weakly closed subspace 
in the space of all bounded linear
operators between $\H_1$ and $\H_2$.}
\endroster

Recall, that the weak topology on the space of bounded linear operators
$f:\H_1\to \H_2$ is given by the family of seminorms 
$$p_{\phi,x}(f)\ =\ |\<\phi,f(x)\>|,\quad \text{where}\quad \phi\in \H_2, 
\quad x\in\H_1.$$

In particular, for any object $\H\in\ob(\ca)$ of a von Neumann category
the set of endomorphisms $\Hom_{\ca}(\H,\H)$ ({\it i.e. the commutant})
is a von Neumann algebra.

Note also that the same notion of von Neumann category will be obtained 
if we require that the set of morphisms $\Hom_{\ca}(\H_1,\H_2)$ is strongly 
(instead if weakly) closed in the space of all bounded 
linear operators between $\H_1$ and $\H_2$. This follows easily from
von Neumann density theorem (cf. \cite{Di}, part I, chapter 3, \S 4)
applied to the von Neumann algebra 
$\Hom_{\ca}(\H_1\oplus\H_2,\H_1\oplus\H_2)$.

Observe, that in a von Neumann category $\ca$ the notions of isomorphism
and unitary equivalence (i.e. $\ca$-isomorphism, preserving the scalar 
products) coincide. Namely, if $f:\H_1\to \H_2$ is an $\ca$-isomorphism then
we can find $T\in\Hom_{\ca}(\H_1,\H_2)$ such that $T^2=f^\ast f$, $T^\ast=T$,
$T>0$ cf. \cite{Di}, part I, chapter 1, \S 2. 
Then the map $g=fT^{-1}:\H_1\to \H_2$ is a unitary equivalence.

\subheading{2.2. Finite objects}
We will say that an object $\H\in\ob(\ca)$ of a von Neumann category 
{\it is finite} if the von Neumann algebra of its endomorphisms
$\Hom_{\ca}(\H,\H)$ is finite. Recall that finiteness of the von Neumann 
algebra $\Hom_{\ca}(\H,\H)$ is equivalent to the fact that {\it any closed 
$\ca$-submodule $\H_1\subset \H$ which is isomorphic to $\H$ in 
$\ca$, coincides with $\H$.} (Cf. \cite{Di}, part III, chapter 8, \S 1.)

If $\ca$ is a von Neumann category, and if $\H\in\ob(\ca)$ is finite 
then any submodule of $\H$ in $\ca$ is also finite. The direct sum of finite
objects is also finite. Thus, all finite objects of any von Neumann 
category form a von Neumann subcategory.

A von Neumann category will be called {\it finite} if all its objects are
finite.

A von Neumann category will be called {\it semi-finite} if any its object 
$\H$ can be represented as a union of an increasing sequence of closed 
finite subobjects $\H_n\subset \H$.

\proclaim{2.3. Lemma} Let $\ca$ be a von Neumann category and let 
$\phi:\H_1\to \H_2$ be an injective morphism of $\ca$ having dense image. 
Then the Hilbert representations $\H_1$ and $\H_2$ are
isomorphic in $\ca$.\endproclaim

\demo{Proof} Consider $f^\ast f: \H_1\to \H_1$; because of property (ii)  
of section (i) we know that
$f^\ast f$ belongs to $\Hom_{\ca}(H_1, H_1)$. Thus $f^\ast f$ is a positive 
element of von Neumann algebra $\Hom_{\ca}(H_1, H_1)$, and so there exists 
$T\in \Hom_{\ca}(H_1, H_1)$ such that
$T^\ast = T$, $T>0$ and $T^2 = f^\ast f$ (cf. \cite{Di}).

Consider the operator $fT^{-1}: TH_1 \to H_2$. It is densely defined and 
bounded. So it can be uniquely extended to a bounded operator 
$fT^{-1}: H_1\to H_2$. Now, we claim that the constructed
operator $fT^{-1}: H_1\to H_2$ belongs to $\Hom_{\ca}(H_1, H_2)$. To show 
this, define
$$A_\epsilon = f\circ \int_\epsilon^\infty \mu_{\epsilon}(\lambda)dE_\lambda$$
where $T=\int_0^\infty \lambda dE_\lambda$ and the real valued continuous
function $\mu_{\epsilon}(\lambda)$ is given by 
$\mu_{\epsilon}(\lambda)=\lambda^{-1}$ 
for $\lambda\ge\epsilon$ and $\mu_{\epsilon}(\lambda)=\epsilon^{-1}$
for $\lambda\le\epsilon$. Then we have:
$A_\epsilon$ belongs to $\Hom_{\ca}(H_1,H_2)$ and $A_\epsilon$ converges to 
$f\circ T^{-1}$
weakly as $\epsilon \to 0$. \qed\enddemo

\proclaim{2.4. Proposition} Any finite object $\H\in\ob(\ca)$ of a 
von Neumann category $\ca$ has the following property.
Suppose that $\phi:\H_1\to \H$ is an injective
morphism in $\ca$ such that its image $\im \phi$ is dense in $\H$. Then 
for any nonzero Hilbert representation $\H_1'\in\ob(\ca)$ and any 
injective morphism $\psi: \H_1'\to \H$ 
in $\ca$, the intersection $\im\phi\cap \im\psi$ is nonzero.
\endproclaim

The role of this property will become clear in the next section where we
will study the torsion objects of the extended category.

\demo{Proof} Suppose that $\ca$ is a von Neumann
category and $\H_2\in\ob(\ca)$ is finite. Let $\phi:\H_1\to\H_2$ be an 
injective morphism of $\ca$ with dense image and let $\psi:\H_1'\to\H_2$ be
another injective morphism of $\ca$ such that $\H_1'$ is nonzero and the 
intersection of the images of $\phi$ and $\psi$ is 0. 

First, we want to find a smaller
nontrivial $\ca$-submodule $\H_1''\subset\H_1'$ such that the image
$\psi(\H_1'')$ is closed. To do so, consider the spectral decomposition
$$\psi^\ast\psi = \int_0^\infty \lambda dE_\lambda$$
and define $\H_1''\ =\ (E_\epsilon\H_1')^\perp $ for some small $\epsilon>0$.
This $\epsilon$ must be chosen so small that $\H_1''$ is nonzero.
Then for any $x\in\H_1''$ we have $|\psi(x)|\ge \sqrt{\epsilon}\cdot|x|$.
Thus, the image of $\psi|_{\H_1''}$ is closed.

Let $\H_2'\subset \H_2$ denote $\psi(\H_1'')^\perp$ and let 
$\pi:\H_2\to\H_2'$ denote the
orthogonal projection. The composite $\pi\circ \phi:\H_1\to \H_2'$ is
injective with dense image. Thus $\H_1$ and $\H_2'$ are isomorphic in $\ca$
by Proposition 2.4. Similarly, applying Proposition 2.4 again we see that
$\H_1$ and $\H_2$ are isomorphic. Thus we obtain that there exists an
isomorphism $\H_2\to\H_2'$ in $\ca$ which contradicts finiteness of $\H_2$.
\qed
\enddemo

\subheading{2.5. Constructing von Neumann categories} We will describe
here two completion constructions which allow to construct 
many interesting examples of von Neumann categories.

\subheading{\it (a) Completing the set of morphisms} Suppose that $\ca$
is a category whose objects have the structure of Hilbert representations
of a $\ast$-algebra $\A$ and such that the morphisms of $\ca$ are represented
by bounded linear maps, i.e. 
$$\Hom_{\ca}(\H_1,\H_2)\subset \L(\H_1,\H_2).$$
We will assume that $\ca$ is {\it self-adjoint}, i.e. condition (ii) of 
section 1.2 is satisfied.

Consider category $\tilde\ca$ which has the same objects as $\ca$ and 
for any two objects $\H_1,\H_2\in\ob(\tilde\ca)$ the set of morphisms
$\Hom_{\tilde \ca}(\H_1,\H_2)$ in $\tilde\ca$ is defines as the closure
of $\Hom_{\ca}(\H_1,\H_2)$ with respect to the weak operator topology.

In order to check that this defines a category we have to show that if
$f:\H_1\to\H_2$ is the limit (in the weak operator topology) of a net
$f_\lambda\in\Hom_{\ca}(\H_1,\H_2)$ and if $g:\H_2\to\H_3$ is the limit
of a net $g_\mu\in\Hom_{\ca}(\H_2,\H_3)$ then the bounded linear map
$g\circ f:\H_1\to\H_3$ belongs the weak closure of $\Hom_{\ca}(\H_1,\H_3)$
in $\L(\H_1,\H_3)$, which we denote $C$. 

We know that the composition of bounded linear maps is {\it not}
continuous considered as a function of two variables (cf.
\cite{Di}, part I, chapter 3); however, it is continuous with respect to
each variable. Thus, we obtain that for any $\mu$ 
$$\lim_\lambda \ g_\mu\circ f_\lambda = g_\mu\circ f$$ 
and so $g_\mu\circ f\in C$
for any $\mu$. Now, using convergence with respect to the other
variable, we obtain that $g_\mu\circ f$ converges to $g\circ f$ and so
$g\circ f$ belongs to $C$.

\subheading{\it (b) Completion by projections}
We will say that an additive category $\ca$ of Hilbert representations of a 
$\ast$-algebra $\A$ is an {\it almost von Neumann category} if it satisfies 
condition (ii) of section 1.2 and condition (v) of section 2.1. 

Our aim in this subsection is to show how
any almost von Neumann category can be canonically completed to a von Neumann
category. 

Let $\ca$ be an almost von Neumann category. Consider a larger category
$\tilde\ca$ of $\ast$-representations of $\A$ whose objects are the
Hilbert representations of the form $e\H$, where $\H\in\ob(\ca)$ and 
$e\in\Hom_{\ca}(\H,\H)$ is a projector $e^\ast=e$, $e^2=e$. 
If $e_1\H_1$ and $e_2\H_2$ are two objects of $\tilde\ca$ then 
{\it a morphism} 
$$f:e_1\H_1\to e_2\H_2$$ 
is defined as a bounded linear map $f$ between the Hilbert spaces 
$e_1\H_1$ and $e_2\H_2$ having the form $f=e_2ge_1$ for some 
$g\in\Hom_{\ca}(\H_1,\H_2)$. Note, that two different maps 
$g_1, g_2\in\Hom_{\ca}(\H_1,\H_2)$ may determine the same map 
$f=e_2g_1e_1=e_2g_2e_1$; in this case they are considered as defining 
the same morphism in $\tilde\ca$.

It is clear that $\tilde\ca$ is a von Neumann category. For example,
the adjoint of a morphism $f=e_2ge_1$ is $f^\ast=e_1g^\ast e_2$; thus 
condition (ii) is satisfied. The kernel of a morphism 
$f:e_1\H_1\to e_2\H_2$, where $f=e_2ge_1$, is $e\H_1$, with 
$e$ denoting $e_1\wedge e_1'$. Here $e_1'$ is the projection onto the 
kernel of $g$ and the wedge denotes the greatest lower bound of $e_1$
and $e_1'$. Note that both $e_1'$ and $e$ belong to $\Hom_{\ca}(\H_1,\H_1)$
by \cite{Di}, part I, chapter 1, \S 2, and by \cite{SZ}, 3.7.

There is the obvious embedding $\ca\to \tilde\ca$ which sends $\H$ to
$1\cdot\H$. Thus {\it we may identify $\ca$ with a full cabcategory of}
$\tilde\ca$.

\subheading{2.6. Examples of von Neumann categories} 

{\bf Example 1.} Let $\A$ be an arbitrary $\ast$-algebra and let $\ca$ be 
the category of all Hilbert representations of $\A$. The objects of this 
category are all
Hilbert representations of $\A$ and the set of morphisms between two given 
representations comprise all bounded linear maps commuting with the $\A$
action. Then $\ca$ is a von Neumann category. It is actually the largest
possible von Neumann category of representations of $\A$. One may also
consider the full subcategory of $\ca$ generated by finite objects; it is
again a von Neumann category.

{\bf Example 2}. The construction which we describe here is a source
of many important examples. Some of them will be discussed below.

Let $\ca$ be a von Neumann category and let $\B$ be a
von Neumann algebra acting on a separable Hilbert space $\H$. 
Denote by $\B'$ the commutant of $\H$. 

Starting from these data we will first describe a category of Hilbert 
representations of the $\ast$-algebra $\A\otimes\B'$ which we denote 
$\cabp$. 

The objects of category $\cabp$ are the {\it Hilbert (completed) tensor
products} of the form $\H_1\otimes\H$ where $\H_1\in\ob(\ca)$
with the obvious action of the algebra $\A\otimes\B'$ (here, abusing the  
notations, the same sign 
$\otimes$ denotes the {\it algebraic} tensor product).

Given two objects $\H_1\otimes \H$ and $\H_2\otimes \H$ of $\cabp$,
the set of morphisms 
$$\Hom_{\cabp}(\H_1\otimes \H, \H_2\otimes \H)\tag2-1$$
is defined as the set of all bounded linear maps 
$\H_1\otimes \H\to\H_2\otimes \H$ of the form
$$\sum_if_i\otimes g_i\tag2-2$$
where the sum is finite and $f_i\in\Hom_{\ca}(\H_1,\H_2)$ and $g_i\in\B$.

The obtained category $\cabp$ is self-adjoint (i.e. it satisfies condition
(ii) of section 1.2) but {\it it may be not a von Neumann category} since it 
may not satisfy (i) of section 1.2 and (v) of section 2.1.

To make a von Neumann category out of $\cabp$ we will apply two completion
constructions described in 2.5. Namely, we will first apply the completion
construction of 2.5.(a) completing the set of morphisms; then we will apply
the construction of 2.5.(b) completing by projectors. The resulting category
will be denoted 
$$\ca\otimes \H.\tag2-3$$
It is a von Neumann category of representations of $\A\otimes \B'$.

{\bf Example 3}. Here we will consider a special case of the previous 
example. We emphasize it because it appears in many applications. It is actually
simpler than the general case since here we do not need to perform 
{\it completing the set of morphisms}, section 2.5(a).

Let $\B$ be a von Neumann algebra acting on a Hilbert 
space $\H$. Denote by $\B'$ be the commutant of $\B$. Consider example 2 
above with $\A=\C$ and with $\ca$ being the category of finite 
dimensional Hilbert spaces. Then the category $\ca\otimes\H$ is a von Neumann
category of Hilbert representations of $\B'$. 

Objects of category $\ca\otimes\H$  are in 
one-to-one correspondence with projections $e\in M(n)\otimes \B$, 
$e^\ast=e, e^2=e$, for some $n$, 
where $M(n)$ is the $n\times n$-matrix algebra.
For each projection $e$ the corresponding Hilbert representation of
$\B'$ is $e(\C^n\otimes\H)$. If $e_1$ and $e_2$ are two projections, 
then the set of
morphisms 
$\Hom_{\ca\otimes\H}(e_1(\C^{n_1}\otimes\H),e_2(\C^{n_2}\otimes\H))$ 
is the set of 
all bounded linear maps of the form
$e_2be_1:e_1(\C^{n_1}\otimes\H)\to e_2(\C^{n_2}\otimes\H)$ 
where $b$ is given by an $n_1\times n_2$-matrix with entries in $\B$. 

Studying this von Neumann category $\ca\otimes\H$ is in some 
sense equivalent to
studying the original von Neumann algebra $\B$. In particular, we see that
this von Neumann category $\ca\otimes\H$ is finite if and only if the original von 
Neumann algebra $\B$ is finite.

{\bf Example 4}. The following example corresponds to the extended category 
constructed in \cite{Fa, Fa1}.

Let $\B$ be a finite von Neumann algebra supplied with a finite, normal and
faithful trace $\tau:\B\to\C$. Consider the completion $\ell^2_{\tau}(\B)$
of $\B$ with respect to the scalar product $\<a,b\>=\tau(ab^\ast)$. Then we 
have $\B$ acting on the Hilbert space $\H = \ell^2_{\tau}(\B)$ from the left 
and from the right, and these actions commute. In fact, to make it 
compatible
with our previous constructions, we may think of the left action of $\B$
on $\H$ and identify the commutant $\B'$ with $\B^\bullet$, the opposite 
algebra of $\B$.

Now, if $\ca$ denotes the von Neumann category of finite dimensional 
Euclidean spaces (here $\A=\C$), then we may consider the von Neumann category
$\ca\otimes\H$ as explained in Example 2. It is a von Neumann category of
$\ast$-representations of the von Neumann algebra $\B^\bullet$.

{\bf Example 5}. Here we consider the special case of example 4 which is most
important for topological applications. 

Let $\pi$ be a discrete group
and let $\H = \ell^2(\pi)$ denote the Hilbert completion of the group ring
$\C[\pi]$. The group ring $\C[\pi]$ acts on $\H$ from both sides. 
We denote by $\B={\Cal N}(\pi)$ the von Neumann algebra of $\pi$
which consists of all bounded linear maps $T:\H\to\H$ commuting with the
left action of $\C[\pi]$. The commutant of ${\Cal N}(\pi)$ can be identified 
with the opposite algebra ${\Cal N}(\pi)^\bullet$.

As above, let $\ca$ denote the von Neumann category of finite dimensional 
Euclidean spaces, then we may form the von Neumann category
$\ca\otimes\H$ as in Example 2. It is a von Neumann category of
$\ast$-representations of ${\Cal N}(\pi)^\bullet$.

{\bf Example 6}. A slightly more general example can be obtained as follows.
Let $V$ be a Hilbert representation of $\pi$ (possibly, finite dimensional)
and let $\H = V\otimes \ell^2(\pi)$ (the completed tensor product). Then 
we have two actions of $\pi$ on $\H$: it acts diagonally from the left
and it acts from the right through the second factor $\ell^2(\pi)$. 
These two actions commute. The right
action clearly extends to an action of the von Neumann algebra 
${\Cal N}(\pi)^\bullet$. The left action may also be extended (by taking the
weak closure) to action of a von Neumann algebra, say $\B$, containing
the group ring $\C[\pi]$. Now we may apply construction of example 3;
the obtained von Neumann category will be a category of 
$\ast$-representations of the algebra ${\Cal N}(\pi)^\bullet$.

{\bf Example 7.}{\it Fields of Hilbert spaces.} This example is important 
for studying families.

Let $Z$ be a locally compact Hausdorff space and let $\mu$ be a 
positive Radon measure on $Z$. Let $\A$ denote the algebra 
$L^\infty_{\C}(Z,\mu)$ (the space of essentially bounded $\mu$-measurable
complex valued functions on $Z$, in which two functions equal locally
almost everywhere, are identical). The involution on $\A$ is given by the
complex conjugation. 

We will construct a category of Hilbert representations
of $\A$ as follows. The objects of $\ca$ are in one-to-one correspondence
with the $\mu$-measurable fields of 
finite-dimensional Hilbert spaces $\xi\to\H(\xi)$ over $(Z,\mu)$, 
cf. \cite{Di}, part II, chapter 1. For any such measurable field of Hilbert
spaces, the corresponding Hilbert space is the direct integral
$$\H = \int^{\oplus}\H(\xi)d\mu(\xi)\tag2-4$$
defined as in \cite{Di}, part II, chapter 1. The algebra $\A$ acts on the 
Hilbert space (2-4) by pointwise multiplication. 

Suppose that we have two $\mu$-measurable finite-dimensional fields of
Hilbert spaces $\xi\to\H(\xi)$ and $\xi\to\H'(\xi)$ over $Z$. Then we 
have two corresponding Hilbert spaces, $\H$ and $\H'$, given as direct
integrals (2-4). We define the set of morphisms $\Hom_{\ca}(\H,\H')$ as the
set of all bounded linear maps $\H\to\H'$ given by {\it decomposable linear
maps} 
$$T = \int^\oplus T(\xi)d\mu(\xi),\tag2-5$$
where $T(\xi)$ is {\it an essentially bounded measurable field of linear 
maps}
$T(\xi):\H(\xi)\to\H'(\xi)$, cf. \cite{Di}, part II, chapter 2.  

The kernel of any decomposable linear map as above can be represented as the
direct integral of a finite-dimensional field of Hilbert spaces and so
the condition (i) of section 1.2 is satisfied. The condition (ii) of section
1.2 is also satisfied since the adjoint of the map $T$ given by (2-5) is
$$T^\ast=\int^\oplus T(\xi)^\ast d\mu(\xi)$$
(by \cite{Di}, part II, chapter 2, \S 3, Proposition 3).
The condition (v) of section 2.1 is satisfied as follows from Theorem 1
of \cite{Di}, part II, chapter 2, \S 5.

Thus, we obtain a von Neumann category of Hilbert representations of 
the algebra $\A = L^\infty_{\C}(Z,\mu)$. 

This category is finite.

Note that this category $\ca$ depends only on the {\it class of the measure
$\mu$}.

The category $\ca$ of this example contains as a full subcategory the 
category $\ca'$ which is obtained by the construction of example 3, 
applied to the Hilbert space $\H=L^2(Z,\mu)$ and the von Neumann algebra
$\B=L^\infty_{\C}(Z,\mu)$ acting on it. Then the commutant $\B'$ is
again $L^\infty_{\C}(Z,\mu)$. In fact, the objects
of the smaller category $\ca'$ can be considered as the measurable fields 
of finite dimensional Hilbert spaces $\xi\to\H(\xi)$ over $Z$ such that 
the dimension of $\H(\xi)$ is essentially bounded.

{\bf Example 8.} We can combine the above examples as follows.
For any von Neumann category $\ca$ and for any
locally compact Hausdorff space $Z$ with a positive Radon measure $\mu$, 
we have a von Neumann category $\ca\otimes L^2(Z,\mu)$ constructed as 
in example 2. It is a category of Hilbert
representations of the algebra $\A\otimes L^\infty(Z,\mu)$. The objects of 
this category may be thought of as {\it "measurable fields of objects of 
the category $\ca$ over the space $Z$, having bounded size",} 
generalizing the previous example 7. 

{\bf Example 9}. {\it Filtrations.} Consider the algebra $\A$ which has 
infinite
number of generators
$e_i$, where $i=1,2,\dots,$ satisfying the relations
$$e_ie_j = e_je_i = e_{\min\{i,j\}}\qquad\text{for all}
\quad i,j = 1,2,3,\dots\tag2-6$$
We define involution on $\A$ by $e_i^\ast = e_i$. 

A Hilbertian representation of $\A$ consists of a Hilbert space $\H$ and
a filtration
$$\H_1\subset\H_2\subset\H_3\subset\dots\subset\H\tag2-7$$
by closed subspaces $\H_i = e_i\H$. 

We denote by $\ca$ the category whose objects are Hilbert representations 
of $\A$ such that all subspaces $\H_i$ are finite dimensional and 
$\H = \cup \H_i$. If $\H$ and
$\H'$ are two such representations, a morphism $\H\to \H'$
in $\ca$ is defined as a bounded linear map commuting with the action of
the algebra $\A$. 

Observe that if $\H$ is an object of $\ca$, then there is the dual 
decreasing filtration
$$\H\supset\H_1^{\perp}\supset\H_2^{\perp}\supset\dots,\tag2-8$$
where $\H_i^{\perp} = (1-e_i)\H$ is the orthogonal complement of $\H_i$.
A bounded linear map $f:\H\to\H'$ is a morphism of $\ca$ if and only if it
preserves the filtration $f(\H_i)\subset \H_i'$ and the dual filtration
$f(\H_i^{\perp})\subset\H_i^{\prime\perp}$. 

It is easy to check that $\ca$ is {\it a finite von Neumann category}.

{\it Remark}. If $\ca$ is a von Neumann algebra and if $\B\subset\A$ 
is a $\ast$-sublagebra then $\ca$ can be obviously viewed as a von 
Neumann category of Hilbertian representations of $\B$.

More generally, if $\phi:\B\to\A$ is a morphism of algebras with involution
then any von Neumann category $\ca$ of Hilbert representations of $\A$ can
be viewed (via $\phi$) as a von Neumann category of Hilbert representations 
of $\B$.

Any full additive subcategory of a von Neumann category 
is again a von Neumann category. Using this remark one may construct other 
examples of von Neumann categories starting from the examples 
mentioned above.

We refer also to \cite{GLR} for other examples.

\subheading{2.7. Traces on von Neumann categories} Let $\ca$ be a 
von Neumann category. 

{\bf Definition.} {\it A trace on $\ca$ is a function, denoted $\tr$, 
which assigns to each object $\H\in\ob(\ca)$ a finite, non-negative trace 
$$\tr_{\H}: \Hom_{\ca}(\H,\H)\to \C$$
on the von Neumann algebra $ \Hom_{\ca}(\H,\H)$; in other words $\tr_{\H}$ assumes
(finite) values in
$\C$, $\tr_{\H}(a)$ is non-negative on non-negative elements $a$ of $ \Hom_{\ca}(\H,\H)$,
and $\tr_{\H}$ is traceful, i.e. 
$\tr_{\H}(ab)=\tr_{\H}(ba)$, for $a, b \in \Hom_{\ca}(\H,\H)$. 

It is also assumed that for any pair of representations
$\H_1$ and $\H_2$ the corresponding traces $\tr_{\H_1}$, $\tr_{\H_2}$ and 
$\tr_{\H_1\oplus\H_2}$ 
are related as follows: if 
$f\in\Hom_{\ca}(\H_1\oplus \H_2, \H_1\oplus \H_2)$ is given by a $2\times 2$
matrix $(f_{ij})$, where $f_{ij}:\H_i\to\H_j$, $i,j=1,2$, then}
$$\tr_{\H_1\oplus\H_2}(f)\ =\ \tr_{\H_1}(f_{11}) + \tr_{\H_2}(f_{22}).
\tag2-9$$

We will say that a trace $\tr$ on a von Neumann category is {\it normal (or faithful)} 
iff for each non-zero $\H\in\ob(\ca)$ the trace $\tr_\H$ on the von Neumann algebra
$\Hom_{\ca}(\H,\H)$ is normal (faithful).

I am going to show elsewhere that non-normal traces produce extremally interesting numerical
invariants of projective and torsion objects of extended abelian categories, which have
important geometric applications.
However in this paper I will consider only normal traces.

{\bf Examples of traces.} Consider the von Neumann category $\ca\otimes \H$ 
of Example 3 of section 2.6. Then any finite trace on the 
von Neumann algebra $\B$ defines a trace on the corresponding von 
Neumann category $\ca\otimes \H$. Namely, given a morphism
$e_2be_1: e_1(\C^n\otimes\H)\to e_1(\C^n\otimes\H)$, where 
$e_1, e_2, b\in M(n)\otimes \B$ and $e_1$ and $e_2$ are projectors
(cf. section 2.6, example 3) we may consider it as an $n\times n$-matrix
$(b_{ij})$ with entries in $\B$ and so we may define
$$\tr(e_2be_1) = \sum_{i=1}^n \tr(b_{ii}).\tag2-10$$

Examples 4, 5, 6 of section 2.6 are special cases of example 3 and so 
using the previous remark we know the traces in these von Neumann categories.

It seems that von Neumann categories of examples 7 and 9 are too 
large to have traces. But some of their subcategories clearly have traces.

For instance, given a locally compact Hausdorff space with a positive Radon
measure $\mu$ one may consider the category of measurable fields of finite
dimensional Hilbert spaces $\H(\xi)$ over $Z$ such that the dimensions of
$\H(\xi)$ are essentially bounded. This category has traces. 
Let $\nu$ be a positive measure on
$Z$ which is absolutely continuous with respect to $\mu$ and such that 
$\nu(Z)<\infty$. Then $\nu$ determines the following trace on this von 
Neumann category
$$\tr_{\nu}(T) = \int_Z \Tr(T(\xi))d\nu\tag2-11$$
where $T$ is a morphism given by formula (2-5) and $\Tr$ denotes the usual
finite dimensional trace.

Let us now describe a von Neumann category which is a full subcategory 
of the category of example 9, which admits a trace. 
Fix a sequence of non-negative real numbers $\mu=(\mu_n)$, such that 
the series $\sum \mu_n$ converges.
Let $\A$ denote the algebra with involution as in example 9. We will 
consider the category $\ca(\mu)$ of all Hilbert
representations $\H$ of $\A$ having the following properties:
\roster
\item the subspaces $e_i\H = \H_i$ are all finite dimensional;
\item $\cup \H_i = \H$;
\item the series $\sum \mu_nd_n$ converges; here $d_n$ denotes the dimension
of the space $\H_{n-1}^\perp\cap \H_n$.
\endroster
Morphisms of $\ca(\mu)$ are all bounded linear maps commuting with $\A$. 

Now we may describe a trace in $\ca(\mu)$. Given an object 
$\H\in\ob(\ca(\mu))$, and a morphism $f:\H\to\H$, define
$$\tr(f) = \sum_{n=1}^\infty \mu_n \Tr(f_n)\tag2-12$$
where $f_n$ is one of the maps between finite dimensional spaces,
$f_n: \H_{n-1}^\perp\cap \H_n \to \H_{n-1}^\perp\cap \H_n$, determined by 
$f$; in other words $f_n = (e_n-e_{n-1})f(e_n-e_{n-1})$. 
The symbol $\Tr$ in the
RHS of (2-12) is the usual finite dimensional trace. Note that the series 
in (2-12) converges for any bounded $f$. It is easy to check that (2-12)
is a trace on the category $\ca(\mu)$.

Note that the notion of trace on a category is also helpful while studying
representations of finite dimensional algebras. If $\A$ is a finite 
dimensional algebra then the category $\ca$ of all finite dimensional
representations of $\A$ is a finite von Neumann category. The traces on
this category can be described as follows. Let $\H_1$, $\H_2$,$\dots, \H_n$
be all different irreducible representations of $\A$. Fix positive numbers
$\mu_1,\dots,\mu_n$. For any $\H\in\ob(\ca)$ and an $\A$-homomorphism 
$f:\H\to\H$ define
$$\tr_{\H}(f) = \sum_{i=1}^n\mu_i\tr(f_i)$$
where $f_i$ denotes the linear map $\Hom_\A(\H_i,\H)\to\Hom_\A(\H_i,\H)$
induced by $f$ and the symbol $\tr$ on the right denotes the usual finite
dimensional trace (which is up to normalization the only trace in the category
of finite dimensional vector spaces).

Note also that it is sometimes more convenient to abandon the condition of
positivity of traces on categories and consider more general traces than
allowed by Definition 2.7.

\subheading{2.8. Von Neumann dimension} Given a trace $\tr$ on a 
finite von Neumann category, one may
measure the sizes of objects of $\ca$ by {\it their von Neumann dimension}
$$\dim \H = \dim_{\tr}(\H) =  \tr_\H(\id_\H)\tag2-13$$
with respect to the chosen trace. The dimensions of non-zero objects of 
$\ca$ are strictly positive if and only if the trace $\tr$ is faithful.

A von Neumann category admitting a faithful trace is necessarily finite.

\proclaim{2.9. Lemma} For any pair of morphisms $f:\H_1\to \H_2$ and
$g:\H_2\to\H_1$ of a von Neumann category $\ca$ supplied 
with a trace $\tr$ holds 
$$\tr_{\H_2}(fg) = \tr_{\H_1}(gf).\tag2-14$$\endproclaim
\demo{Proof} 
$$
\aligned
\tr_{\H_1}(gf) =& \tr_{\H_1\oplus\H_2}(\bmatrix gf & 0\\ 0 & 0\endbmatrix)\\
=& \tr_{\H_1\oplus\H_2}(\bmatrix 0 & 0\\ g & 0\endbmatrix\circ \bmatrix 0 & 
f\\ 0 & 0\endbmatrix)\\
=& \tr_{\H_1\oplus\H_2}(\bmatrix 0 & f\\ 0 & 0\endbmatrix\circ \bmatrix 0 & 0\\
g & 0\endbmatrix)\\
=& \tr_{\H_1\oplus\H_2}(\bmatrix 0 & 0\\ 0 & fg\endbmatrix) \\
=& \tr_{\H_2}(fg)
\endaligned
$$
\enddemo

\proclaim{2.10. Corollary} If $\H_1, \H_2 \in\ob(\ca)$ are isomorphic in $\ca$ 
then their von Neumann dimensions are equal
$\dim(\H_1) = \dim(\H_2)$
\endproclaim
\demo{Proof} Let $f:\H_1\to \H_2$ and $g:\H_2\to\H_1$ be mutually inverse
isomorphisms. Then 
$$
\dim \H_1 = \tr_{\H_1}(\id_{\H_1}) = \tr_{\H_1}(gf) =  
\tr_{\H_2}(fg) = \tr_{\H_2}(\id_{\H_2}) = \dim \H_2
$$
\enddemo

It is well known that any finite von Neumann algebra admits a finite normal faithful
trace. It is probably not true that any finite von Neumann category admits a faithful
trace.
 
\proclaim{2.11. Question} Which conditions guarantee existence
of faithful traces on a given finite von Neumann category?
\endproclaim

\heading{\bf \S 3. Torsion objects of the extended category}\endheading

In this section we will further study the construction of the extended
abelian category $\eca$ of section 1. 

We will here assume everywhere in this section that {\it the original 
category $\ca$ is a finite von Neumann category}.

\subheading{3.1. Definition} An object $\X=(\alpha:A'\to A)$ of the extended
category $\eca$ is called {\it torsion} iff the image of $\alpha$ is dense 
in $A$.

Equivalently, an object $\X$ of $\eca$ is torsion if and only if it admits 
no nontrivial homomorphisms in projective objects of $\eca$.

We will denote by $\tca$ the full subcategory of $\eca$ generated by all
torsion objects. $\tca$ is called {\it the torsion subcategory of $\eca$}.

Note, that in general (without assuming that $\ca$ is a finite von Neumann
category) it is possible that a projective object is subobject of a 
torsion object.  Our desire to avoid this situation explains the 
assumption on $\ca$.

\proclaim{3.2 Proposition} Given an exact sequence
$$0\to \X'\to \X\to \X''\to 0\tag3-1$$
of objects and morphisms of $\eca$, where $\ca$ is a finite von Neumann 
category, the middle object $\X$ is torsion if and
only if both $\X'$ and $\X''$ are torsion.\endproclaim
\demo{Proof} The only nontrivial part of the proof consists in showing that
a subobject of a torsion object is necessarily torsion. As we will see, 
this follows from Proposition 2.4. 

Suppose that the diagram
$$
\CD
(\alpha:A'@>>> A)\\
 @.          @VVfV\\
(\beta:B'@>>>B)
\endCD
$$
represents a monomorphism in $\eca$ and $(\beta:B'\to B)$ is torsion.
Then from Proposition 1.6 it follows that 
$$\alpha(A')\supset f^{-1}(\beta(B').\tag3-2$$ 

Suppose that $(\alpha:A'\to A)$ is 
not torsion. Then there is a nontrivial closed submodule $C\subset A$ such 
that $\cl(\im(\alpha))\cap C = 0$. We obtain from (3-2) that $f$ maps $C\to B$
monomorphically and that $f(C)\cap \beta(B')=0$. This contradicts Proposition
2.4. $\qed$
\enddemo

\proclaim{3.3. Corollary} The torsion subcategory $\tca$ is an abelian 
subcategory of $\eca$.
\endproclaim

Recall that we assume that $\ca$ is a finite von Neumann category.

\subheading{3.4. Functor of the torsion part} Given an arbitrary object
$\X=(\alpha:A'\to A)$ of $\eca$ one may consider the following torsion
object 
$$\T(\X) =\ (\alpha: A'\to \cl(\im(\alpha)))\tag3-3$$
which is called {\it the torsion part of $\X$.} There is an obvious morphism
$\T(\X)\to\X$ which is a monomorphism and so $\T(\X)$ can be viewed as a 
subobject of $\X$. It is clear that any morphism
$\X\to\Y$ in $\eca$ maps the torsion part of $\X$ into the torsion part of
$\Y$. Thus, $\T$ is naturally defined as a functor
$$\T: \eca \to \ca.\tag3-4$$

Using Proposition 1.6, we easily see that, the following sequence is exact
$$0\to \T(\X) \to \X \to P(\X)\ \to 0.\tag3-5$$
Here $P$ denote the functor of the projective part which was discussed in
1.11. Since $P(\X)$ is projective (by Proposition 1.9) the sequence (3-5)
splits and so we have
$$\X = \T(\X) \oplus P(\X).\tag3-6$$
We obtain: {\it isomorphism type of an object of the extended category $\eca$
is determined by the isomorphism types of its projective and torsion parts.}

Note that representation (3-6) is not canonical.

\subheading{3.5} Our next aim is to apply the construction of S.Novikov and
M.Shubin \cite{NS} which allows to measure the sizes of torsion objects 
of the extended category. 

We start with the following technical statement which will be quite 
important for the sequel.

We will say that two torsion objects $\X=(\alpha:A'\to A)$ and 
$\Y=(\beta:B'\to B)$ of the extended category $\eca$ (with injective $\alpha$
and $\beta$) are {\it strongly isomorphic} is there exist two isomorphisms
$f:A\to B$ and $g:A'\to B'$ in $\ca$ such that $f\circ\alpha = \beta\circ g$.
Note that strongly isomorphic objects are necessarily isomorphic, but the
converse is not true.

\proclaim{3.6. Proposition} Let $\X=(\alpha:A'\to A)$ and 
$\Y=(\beta:B'\to B)$ be
two isomorphic torsion objects of the extended category $\eca$ 
with injective $\alpha$ and $\beta$. Then
it is possible to perform a finite number of excisions on $\X$ and on $\Y$ 
(at most four!) such that
the obtained objects will be strongly isomorphic.
\endproclaim
\demo{Proof} Suppose that the diagram
$$
\CD
(\alpha:A'@>>> A)\\
      @VgVV   @VVfV\\
(\beta:B'@>>>B)
\endCD\tag3-7
$$
represents an isomorphism in $\eca$. From Proposition 1.6 it follows 
(using the fact that the morphism represented by $f$ is a monomorphism)
that $\alpha(A')\supset f^{-1}(\beta(B'))$. In particular, the kernel of
$f$ is contained in the image of $\alpha$. Thus, we may perform an excision
with respect to $\alpha^{-1}(\ker(f))$. As a result we will obtain a 
diagram similar to (3-7), which represents an isomorphism in $\eca$ with
$f$ and $g$ injective. 

Consider now $\im(f)^\perp\subset B$ and the orthogonal projection
$\pi: B\to \im(f)^\perp$. From Proposition 1.6 it follows 
(using the fact that the morphism represented by $f$ is an epimorphism)
that $\im(f)^\perp$ coincides with the image of $\pi\circ\beta$. Thus we may 
perform an excision on $(\pi\circ\beta)^{-1}(\im(f)^\perp)$ such that the 
initial isomorphism
will be represented by a diagram of form (3-7) such that $f$ and $g$ are
injective with dense images.

By Proposition 2.3 all Hilbert representations $A$, $B$, $A'$, $B'$ are
isomorphic in $\ca$ and so we may identify them. So, we may assume that we 
are given an isomorphism in $\eca$ of the form
$$
\CD
(\alpha:A@>>> A)\\
 @VgVV          @VVfV\\
(\beta:A@>>>A)
\endCD \tag3-8
$$
with $\alpha$, $\beta$, $f$ and $g$ being injective with dense images.
Also, we may suppose that $\alpha$ and $\beta$ are self-adjoint and positive.

By Proposition 1.6, the sequence
$$0\to A @>{\alpha\oplus -g}>> A\oplus A @>{(f,\beta)}>> A\to 0$$
is exact in $\ca$ and so it splits. Thus there exist morphisms
$$\sigma, \delta: A\to A$$
such that 
$$\id_A = \sigma \alpha + \delta g. \tag3-9$$
Choose $\epsilon>0$ such that $\epsilon \cdot||\sigma||<1$. 

Now, we may find a splitting
$$A = P \oplus Q$$
such that $\alpha(P) = P$ and $\alpha(Q)\subset Q$ and the norm of the
restriction of $\alpha$ on $Q$ will be less than $\epsilon$. This may be
achieved by choosing a corresponding spectral projector
from the spectral decomposition of $\alpha$. Then from equation (3-9) we find
that $\pi_Q\circ\delta g|_Q:Q\to Q$ is an isomorphism (where $\pi_Q:A\to Q$
denotes the orthogonal projection) and thus $g|_Q:Q\to A$ has closed image.
Hence we obtain that the initial isomorphism will be represented by a
diagram
$$
\CD
(\alpha|_Q: Q@>>> Q)\\
 @V{g|_Q}VV          @VV{f|_Q}V\\
(\beta:A@>>>A)
\endCD\tag3-10
$$
such that the image of $g|_Q$ is closed. Then the image of $f|_Q$ is closed
(here one again uses Proposition 1.6 and the fact that the diagram 
(3-10) represents an isomorphism). As above, we may perform an excision with 
respect to the orthogonal complement to $f(Q)$ in $A$ to obtain a strong
isomorphism. \qed
\enddemo

\subheading{3.7. The spectral density function} We are now ready to 
describe the construction which was used by S.Novikov and M.Shubin \cite{NS}
in order to measure the spectrum near zero. In the language of the present
paper this construction assigns to
a torsion object of the extended category a {\it spectral density function,}
which is real valued function determined up to certain equivalence relation. 
There is also an
important numerical invariant which can be extracted from the spectral 
density function (it is called {\it the Novikov-Shubin number}, cf. 3.9).

Suppose that $\ca$
is {\it a finite von Neumann category supplied with a trace} $\tr$.

Let $\X=(\alpha:A'\to A)$ be a torsion object of the extended category
represented by an injective morphism $\alpha$. Let $T:A'\to A'$ be a
positive self-adjoint operator, $T\in\Hom_{\ca}(A',A')$, such that
$T^2=\alpha^\ast\alpha$. By the spectral theorem, we have representation
$$T = \int_0^\infty\lambda dE_\lambda\tag3-11$$
and we denote by $F(\lambda)$ the von Neumann dimension (with respect to the
chosen trace $\tr$) of the subspace $E_\lambda A'$ determined by the spectral
projector $E_\lambda$, where $\lambda>0$: 
$$F(\lambda) = \dim E_\lambda A'\tag3-12$$
The obtained function $F(\lambda)$
is called {\it the spectral density function of $\X$}.

It is clear that $F(\lambda)$ is non-decreasing and right continuous (since 
the trace $\tr$ is supposed to be normal) function with $F(0)=0$. If the
trace $\tr$ is faithful then $F(\lambda)>0$ for $\lambda>0$. Note, that in
general we do not require the trace $\tr$ to be faithful.

Let us emphasize again that {\it the spectral density function of 
a torsion object $\X$ depends on $\X$ and on the trace $\tr$} although 
our notation does not indicate this fact.

Since we are interested in getting an isomorphism invariant of
$\X$, we will focus our attention only in the behavior of the spectral 
density function {\it near zero}; in fact, performing excisions, we may 
change the spectral density function $F(\lambda)$ arbitrarily away from zero. 

Two spectral density functions $F(\lambda)$ and $G(\lambda)$ are called
{\it dilatationally equivalent} (denoted $F(\lambda)\sim G(\lambda)$)
if there exist constants $C>1$ and 
$\epsilon>0$ such that
$$G(C^{-1}\lambda) \le F(\lambda) \le G(C\lambda)$$
for all $\lambda\in (0,\epsilon)$.

The following statement was proven by M.Gromov and M.Shubin \cite{GS}
in a much more general situation. We include its proof for the sake of
completeness.

\proclaim{3.8. Proposition} Spectral density functions of isomorphic 
torsion objects $\X$ and $\Y$ are dilatationally equivalent.
\endproclaim
\demo{Proof} Because of Proposition 3.6, it is enough to show that the
dilatational equivalence class of the spectral density function does not
change under excisions and under strong isomorphisms.

Excisions clearly change the spectral density function only away from
zero. Thus we need only to consider strong isomorphisms. Suppose that a
commutative diagram
$$
\CD
(\alpha:A'@>>> A)\\
 @VgVV          @VVfV\\
(\beta:B'@>>>B)
\endCD\tag3-14
$$
(with $\alpha$ and $\beta$ injective) represents a strong isomorphism (i.e.
$f$ and $g$ are isomorphisms). Let $F(\lambda)$ and $G(\lambda)$ denote 
the spectral density functions of $\X= (\alpha:A'@>>> A)$ and 
$\Y= (\beta:B'@>>>B)$ correspondingly. Since $f$ and $g$ are isomorphisms, 
there exists a constant $K>1$ such that
$$K^{-1}|x| \le |f(x)| \le K|x|\quad\text{and also}\quad 
K^{-1}|x| \le |g(x)| \le K|x|\tag3-15$$
for all $x$. 

We will use the variational characterization of the spectral density 
function given by M.Gromov and M.Shubin \cite{GS}. Namely, for any
$\lambda>0$ consider the cone
$$C_\lambda (\X) = \{x\in A'; |\alpha(x)|\le \lambda |x|\}.\tag3-16$$
Then 
$$F(\lambda) = \sup \dim L,\tag3-17$$
where $L$ runs over all $\ca$-subspaces $L\subset A'$, $L\in \ob(\ca)$.

Now, using inequalities (3-15) we obtain that if $L\subset C_\lambda(\X)$
then $gL\subset C_{K^2\lambda}(\Y)$. Combined with (3-17), it implies
$G(K^2\lambda)\ge F(\lambda)$. 

Changing the role of $\X$ and $\Y$, one obtains that 
$F(K^2\lambda)\ge G(\lambda)$. \qed
\enddemo

\subheading{3.9. The Novikov - Shubin number} In this subsection we will
assume that the trace $\tr$ is normal.

Let $\X$ be a torsion object
of the extended category $\eca$. Suppose that its spectral density function
$F(\lambda)$ constructed with respect is a given trace $\tr$ on the initial
von Neumann category $\ca$ is {\it positive} for $\lambda>0$. Then one 
defines the following
non-negative real number (which allowed to be $\infty$)

$$\ns(\X) = \lim_{\lambda\to 0^+} \inf \frac{\log F(\lambda)}{\log \lambda}
\ \in\ [0,\infty],\tag3-18$$
which is called {\it the Novikov - Shubin number or the Novikov - Shubin 
invariant of $\X$.} Sometimes we will use the notation
$$\ns(\X) = \ns_{\tr}(\X)\tag3-19$$
in order to emphasize dependence on the trace $\tr$.

From Proposition 3.8 it follows that {\it isomorphic torsion objects 
$\X \simeq \Y$ have equal Novikov - Shubin numbers $\ns(\X) = \ns(\Y)$}. 

In fact, it is more convenient to use an equivalent invariant
$$\c(\X) = \ns(\X)^{-1}\ \in\ [0,\infty],\tag3-20$$
which was introduced in \cite{Fa1}; it is called {\it the capacity of} 
$\X$. The advantage of capacity against the Novikov-Shubin number
consists in the fact that it adequately describes the size of a torsion 
object: larger torsion objects have larger capacity, trivial object 
has zero capacity, there may also exist torsion objects having infinite
capacity.

If the spectral density function $F(\lambda)$ vanishes for some $\lambda>0$ 
then we will define $\c(\X) = 0$. This situation may happen for a nonzero
$\X$ if the trace is not faithful.

The following statement describes the main properties of the Novikov - Shubin
capacity.

\proclaim{3.10. Proposition} For any short exact sequence of torsion
objects 
$$0\to \X' \to \X \to \X'' \to 0$$
holds
$$\max\{\c(\X'),\c(\X'')\}\le \c(\X) \le \c(\X')+\c(\X'')\tag3-21$$
and 
$$ \c(\X) =  \max\{\c(\X'),\c(\X'')\}\tag3-22$$
if the exact sequence above splits. In particular,
$$\c(\X\oplus\X) = \c(\X)\tag3-23$$
for any torsion $\X$.
\endproclaim

Cf. \cite{Fa1}, Proposition 4.9. The arguments of \cite{Fa1} are in fact 
based on results of \cite{LL}.

\proclaim{3.11. Corollary} For a fixed trace $\tr$ on a finite von Neumann
category, and for a given $\nu\in [0,\infty]$, let $\T_\nu(\ca)$ denote the 
full
subcategory of $\tca$ whose objects are all torsion objects $\X$ of $\eca$
with capacity (taken with respect to the trace $\tr$) $\c(\X)\le \nu$. Then
$\T_\nu(\ca)$ is an abelian subcategory.\qed\endproclaim

\subheading{3.12. Example} {\it In a finite von Neumann category there may 
exist non-isomorphic torsion objects which have identical Novikov-Shubin 
numbers with respect to any trace on the initial von Neumann category.}

To show this, consider the von Neumann category $\ca$ of measurable fields of
finite dimensional Hilbert spaces over the circle $S^1$ supplied with the 
Lebesgue measure. The algebra $\A = L^\infty(S^1)$ clearly acts on the 
objects of this category $\ca$. The traces on this 
category are in one-to-one correspondence with the measures on the 
circle $S^1$, which are absolutely continuous with respect to the Lebesgue
measure. 

Let $\H\in\ob(\ca)$ denote the object represented by a constant field of 
one-dimensio\-nal spaces over the circle. The corresponding Hilbert space
of $L^2$ sections is $\H = L^2(S^1)$. For any angle $\theta\in [0,2\pi]$
and for any real $\nu >0$ consider the following torsion object
$\X = \X_{\theta,\nu} = (\alpha :\H\to \H)$, where $\alpha$ is the operator 
of pointwise multiplication by the function 
$\alpha(z) = |z-e^{i\theta}|^\nu$, where $z\in S^1$. 

It follows from the claim in section 7.4 of \cite{Fa1} that $\X$ and
$\X\oplus\X$ are not isomorphic. But they have identical capacities 
with respect to any trace according to Proposition 3.10 above.

\subheading{3.13. Duality for torsion objects} Let $\ca$ be a finite 
von Neumann category, and let $\X=(\alpha:A'\to A)$ be a torsion
object, such that $\alpha$ is injective. We will define the {\it dual object}
$\e(\X)$ as $\e(\X)=(\alpha^\ast:A\to A')$. Here the star denotes the
adjoint morphism. If $\Y=(\beta:B'\to B)$ is another
torsion object with injective $\beta$ then any morphism 
$[f]:\X\to\Y$ represented by a morphism $f:A\to B$ in $\ca$, determines the
dual morphism $\e([f]):\e(\Y)\to\e(\X)$ which is defined as follows.
According to the definition of section 1.3, there exists a morphism
$g:A'\to B'$ in $\ca$ such that $f\circ \alpha = \beta\circ g$; this $g$ is
in fact unique since $\alpha$ and $\beta$ are supposed to be injective.
Then $\e([f])$ is defined as morphism $[g^\ast]$.

The duality $\e$ is a contravariant functor $\e:\T(\ca)\to\T(\ca)$,
and there is natural equivalence $\e\circ \e\simeq \Id$.

We refer to \cite{Fa1}, section 3.8 for a more detailed description of this 
duality in a specific situation.

Observe that {\it any torsion object is isomorphic to its dual}, but this 
isomorphism
is not canonical. In fact, by Proposition 2.3, we know that any torsion
object is isomorphic to $(\alpha:A\to A)$; then using the polar decomposition
we may replace $\alpha$ by a self-adjoint map $\alpha^\ast=\alpha$. In 
this case the dual of $\X=(\alpha:A\to A)$ is identical to $\X$.

\subheading{3.14. Computing the Novikov - Shubin invariants of cohomology}
In many geometric applications torsion objects appear as cohomology
of projective cochain complexes. We need to be able to compute their
spectral density functions and the Novikov - Shubin numbers in terms of the
given chain complex. We are going to make some general remarks
concerning this question, which will be used later in applications.

Suppose that $\ca$ is a finite von Neumann category and $\tr$ is a
fixed trace in $\ca$. Let 
$$C^\ast = (\dots\to C^{i-1}@>{\partial}>> C^{i}@>{\partial}>> 
C^{i+1}@>{\partial}>> \dots)\tag3-24$$
be a cochain complex in $\ca$. By Proposition 1.10 we have the following
formula for the extended $L^2$ cohomology
$$\H^i(C) = (\partial: C^{i-1}\to Z^i)\tag3-25$$
where $Z^i =\ker[\partial:C^i\to C^{i+1}]$. Applying Definition 3.7, we 
obtain the following recipe of computing the spectral density function
of the torsion part of $\H^i(C)$. 
Consider the self-adjoint operator
$$\partial^\ast\partial: C^{i-1} \to C^{i-1}\tag3-26$$
and let
$$G(\lambda) = \tr(\chi_{[0,\lambda]}(\partial^\ast\partial)),\tag3-27$$
where $\chi_{[0,\lambda]}$ denotes the 
characteristic function of the interval $[0,\lambda]$.
Then {\it the von Neumann dimension of the projective part of $\H^i(C)$ 
equals
$$\dim_{\tr}P(\H^i(C)) = G(0)\tag3-28$$
and the spectral density function of the torsion part of $\H^i(C)$
equals} 
$$F^i(\lambda) = G(\lambda^2) - G(0).\tag3-29$$

We shall also point out the relation between the spectral density function
of the Laplacian 
$\Delta^i = \partial\partial^\ast +\partial^\ast\partial :C^i \to C^i$
and the functions $F^i(\lambda)$ which appear above.
Namely, 
$$\tr(\chi_{[0,\lambda]}(\Delta^i)) = F^i(\sqrt{\lambda}) + 
F^{i+1}(\sqrt{\lambda}).\tag3-30$$
We will leave the proof of this formula as an easy exercise.

\heading{\bf \S 4. Divisor determined by a torsion object}\endheading 

This section devoted to studying {\it "families"} or {\it "the abelian
case"}. Namely,
let $Z$ be a locally compact Hausdorff space and let $\mu$ be a positive
Radon measure on $Z$. In this section we will deal with a von 
Neumann category $\ca$ of Hilbert representations of the algebra 
$\A = L^\infty_\C(Z,\mu)$ of essentially bounded complex valued
functions on $Z$. We will see that any  
torsion object $\X$ of the corresponding extended category $\eca$ 
determines a "divisor" - a closed subset $\D(\X)\subset Z$ with additional 
information on it ({\it "multiplicities"}). We compute explicitly the
divisor assuming that the torsion object is given as the cohomology of a 
chain complex of vector bundles. 

We compute also some simple examples which provide an evidence that by 
considering torsion objects of the extended category it is possible to 
capture all relevant information about intersections of analytic varieties.

\subheading{4.1} Let $\ca$ denote the following von Neumann category 
of Hilbert representations of $\ast$-algebra $\A=L^\infty_\C(Z,\mu)$ 
The objects of category $\ca$ are
in one-to-one correspondence with measurable fields of finite dimensional 
Hilbert spaces $\xi\to\H(\xi)$ over $Z$ such that the dimension $\dim\H(\xi)$ 
is essentially bounded. Any such field $\H(\xi)$ determines a Hilbert
representation of $\A$ acting on the direct integral of Hilbert spaces (2-4). 
The morphisms of $\ca$ are bounded linear maps between the corresponding
Hilbert spaces given by decomposable linear maps (2-5), cf. example 7, 
section 2.6 (this category was denoted $\ca'$ then). 

From section 2.7 we know that any measure $\mu'$ on $Z$ which is absolutely
continuous with respect to $\mu$ and such that $\mu'(Z)<\infty$ determines
a trace $\tr_{\mu'}$ on von Neumann category $\ca$, cf. (2-11).

Let $\eca$ denote the extended category of $\ca$ 
and let $\X$ be a torsion object of $\eca$. 
Applying construction of 
section 3.7, we obtain that 
any measure $\mu'$ as in the previous paragraph determines a 
trace $\tr_{\mu'}$ on category $\ca$ and,
the corresponding spectral density function of 
$\X$ with respect to $\tr_{\mu'}$; we denote it by $F_{\mu'}(\lambda)$. 

\subheading{4.2} Given a torsion object $\X$ of $\eca$ we define {\it its
divisor} $\D(\X)\subset Z$ as the set of all points $\xi\in Z$ such that 
for any
neighbourhood $V$ of $\xi$ there exists a measure $\mu'\ll\mu$ 
with support of $\mu'$
contained in $V$ and such that the corresponding spectral
density function $F_{\mu'}(\lambda)$ is {\it strictly positive} 
for all $\lambda>0$.

$\D(\X)$ is a closed set. To see this, one just forms the negation of the 
definition in the previous paragraph: $\xi\notin\D(\X)$ iff there exists a
neighbourhood $V$ of $\xi$ such that for any measure $\mu'\ll\mu$ with 
support in $V$ the corresponding spectral density function 
$F_{\mu'}(\lambda)$
of $\X$ vanishes for all small $\lambda>0$.

Proposition 3.8 implies that 
{\it any two isomorphic torsion objects $\X$ and $\Y$ have identical
divisors} $\D(\X) = \D(\Y)$.

The following statement computes the divisor in many interesting 
situations.

\proclaim{4.3. Proposition} Let $Z$ be a compact Hausdorff 
space supplied with a positive Radon measure $\mu$ such that there are no
nonempty open sets $U\subset Z$ with $\mu(U)=0$. 
Let $\E$ be a Hermitian vector bundle over $Z$ of rank $n$ 
(i.e. a locally trivial vector bundle such that
each fiber is given a Hermitian metric which varies continuously). 
Let $T:\E\to\E$ be a continuous map
which induces a linear self-map on each fiber $\E_\xi$, where $\xi\in Z$. 
For any point
$\xi\in Z$ denote by $d(\xi)$ the determinant of the linear map 
$T(\xi):\E_\xi\to\E_\xi$ of the fiber above $\xi$. Suppose that the zero
set $\{\xi\in Z; d(\xi) = 0\}$ of the obtained continuous 
function $d:Z\to\C$ has measure zero. Then the induced
by $T$ map $T^\ast$ of the spaces of $L^2$ sections 
$\X = (T^\ast: L^2(\E) \to L^2(\E))$
represents a torsion object of the extended category $\eca$
and the divisor $\D(\X)$ of $\X$
coincides with the zero set $\{\xi\in Z; d(\xi) = 0\}$ of the function $d$.
\endproclaim

\demo{Proof} The linear map $T(\xi)$ is invertible $S(\xi)=T^{-1}(\xi)$
for $d(\xi)\ne 0$. We can find a sequence of open sets 
$U_1\supset U_2\supset U_3\supset\dots$ with $\mu(U_j)\to 0$, which
contain the zero set
$\{\xi\in Z; d(\xi) = 0\}$ and such that  all the spaces
$L^2(\E|_{Z-U_n})$ belong to the image of $T^\ast$ and their union
$\cup L^2(\E|_{Z-U_n})$ is dense in $L^2(\E)$. This shows that the map 
$T^\ast$ has dense image and so $\X$ is a torsion object.

Suppose that $\xi\in Z$ is such that $d(\xi)\ne 0$. Then there is a 
neighborhood
$U$ of $\xi$ such that the bundle map $T|_U:\E|_U\to\E|_U$ is invertible. 
Then clearly, for any measure $\mu'\ll\mu$ with $\supp(\mu')\subset U$ we will
have $F_{\mu'}(\lambda)=0$ for small $\lambda>0$. Thus $U\cap\D(\X)=\emptyset$.
This proves that the divisor $\D(\X)$ is contained in the zero set of $d$.

Suppose now that $\xi\in Z$ and $d(\xi)=0$. We want to show that for any
closed 
neighborhood $U$ of $\xi$ the spectral density function $F_{\mu'}(\lambda)$
is strictly positive for all $\lambda>0$, where $\mu' = \mu|_U$. Suppose the 
contrary,
i.e. there exists a closed neighborhood $U$ of $\xi$ and a positive number
$\lambda_U$ such that taking $\mu'=\mu|_{U}$ we have for the corresponding
spectral density function $F_{\mu'}(\lambda) = 0$ for all 
$0< \lambda <\lambda_U$.
Restricting the neighborhood $U$ if necessary, we may additionally assume
that 
$$|d(\xi)| < \lambda_U^n\tag4-1$$
for all $\xi\in U$. Also, we may assume that the bundle $\E$ can be 
trivialized over $U$ and so $T(\xi)$ is given by an $n\times n$-matrix 
$(a_{ij}(\xi))$ of continuous functions of $\xi\in U$. Consider the 
non-negative self-adjoint matrix 
$$(T^\ast T)(\xi) = (c_{ij}(\xi)),\quad\text{where}\quad
c_{ij}(\xi) = \sum_{k=1}^n \overline a_{ik}(\xi)a_{kj}(\xi).\tag4-2$$
We may find measurable functions 
$\lambda_i:U\to \R$, where $i=1,2, \dots, n$, such that the numbers
$$0\le \lambda_1(\xi)\le \lambda_2(\xi)\le \dots \le \lambda_n(\xi)\tag4-3$$
form the set of eigenvalues of the matrix (4-2) for almost all $\xi\in U$.
Applying the definition of spectral density function in section 3.7, we
obtain
$$F_{\mu'}(\lambda) = \sum_{i=1}^n \mu\{\xi\in U;\lambda_i(\xi)\le \lambda^2\}.
\tag4-4$$
Now, our assumption that $F_{\mu'}(\lambda)=0$ for all $\lambda<\lambda_U$
implies that for any $i=1,2,\dots, n$ holds
$$\mu\{\xi\in U; \lambda_i(\xi)>\lambda_U^2\} = \mu(U),\tag4-5$$
and so 
$$\mu(\bigcap_{i=1}^n\{\xi\in U; \lambda_i(\xi)>\lambda_U^2\}) = \mu(U).\tag4-6$$  
However 
$$d(\xi)^2 = \prod_{i=1}^n \lambda_i(\xi) < \lambda_U^{2n}\tag4-7$$
almost everywhere on $U$, which gives a contradiction. \qed
\enddemo

\subheading{4.4. Example} Consider the special case of the situation of the
previous Proposition when $\E$ is the trivial line bundle. Then the bundle
map $T$ is given by a continuous function $f:Z\to \C$ and the corresponding
operator $T^\ast: L^2(Z)\to L^2(Z)$ is multiplication by $f$, which we
will denote $M_f$. Proposition 4.3
states that in this case the divisor of the torsion object 
$\X = \X_f = (M_f: L^2(Z)\to L^2(Z))$ equals to the zero set 
$\{\xi\in Z;f(\xi)=0\}$, provided the measure of the zero set with respect to
$\mu$ vanishes.

Let us take, for instance, $Z=S^1$, and the measure $\mu$ being the Lebesgue
measure. Let $f:S^1\to\C$ be the following function
$f(z)=|z-e^{i\theta}|^\nu$, where $\theta\in [0,2\pi]$ and $\nu>0$. Denote
the corresponding torsion object $\X_{\theta,\nu}$ (it appeared already
in section 3.12.) The divisor
$\D(\X_{\theta,\nu})$ is the one-point set $\{\theta\}$. This shows that
for $\theta\ne\theta'$ the objects $\X_{\theta,\nu}$
and $\X_{\theta',\nu'}$ are not-isomorphic, cf. also \cite{Fa1}. 

Note also that the Novikov - Shubin number $\ns(\X_{\theta,\nu})$ with respect
to the Lebes\-gue measure on $S^1$ equals $\nu^{-1}$, cf. \cite{Fa1}. 
Thus, one may
find the number $\nu$ by means of the Novikov - Shubin number.

\proclaim{4.5. Corollary} In the situation of Proposition 4.3, the condition 
$\D(\X) = \emptyset$ implies $\X=0$.\endproclaim

Indeed, $\D(\X)=0$ implies, by Proposition 4.3, that the determinant function
$d$ nowhere vanishes and thus $T(\xi)$ is invertible for all $\xi$. 
Because $Z$ is compact we obtain that the induced map on the spaces of 
$L^2$-sections $T^\ast:L^2(\E)\to L^2(\E)$ is invertible and so $\X=0$. 

Note that this corollary is not true if $Z$ is only locally compact. 
Consider the 
following example. Let $Z=\R$ with the measure $\mu = dx/x^2$, and let 
$f(x)=e^{-x^2}$. Since $f$ is bounded, it induces a map between the 
spaces of $L^2$-section. Then the torsion object 
$\X_f=(M_f: L^2(Z,\mu)\to L^2(Z,\mu))$ has empty divisor $\D(\X)=\emptyset$,
but $\X\ne 0$.

This suggests the idea to study {\it "the divisor at infinity"} as well.

\proclaim{4.6. Proposition}{\it (Integrality of the Novikov - Shubin 
capacity)} Let 
$Z$ be a compact smooth manifold and let $f:Z\to \R$ be a smooth function 
which is 
transversal with respect to $0\in \R$.
For any positive integer $m>0$ consider the following torsion object
$$\X_{f^m} = (M_{f^m}: L^2(Z)\to L^2(Z)),$$ 
where the $L^2(Z)$ is constructed with respect to a smooth measure on $Z$.
Then the Novikov - Shubin capacity $\c(\X_{f^m})$ with respect to any
smooth measure $\mu$ on $Z$ is integral and equals $m$. Moreover, the spectral
density function of $\X_{f^m}$ is dilatationally equivalent to $\lambda^{1/m}$.
\endproclaim

\demo{Proof} Let $\D = \D(\X_{f^m}) = \{\xi\in Z; f(\xi) = 0\}$ denote 
the divisor of 
$\X$. It is a smooth codimension one submanifold of $Z$. The spectral density 
function
$F_\mu(\lambda)$ of $\X_{f^m}$ is
$$F_\mu(\lambda) = \mu(U_\lambda),\quad\text{where}\quad 
U_\lambda = \{\xi\in Z; |f(\xi)|^m \le \lambda\}.\tag4-8$$
From this formula it clear that it is enough to prove Proposition 4.6 in case
$m=1$; therefore we will assume that $m=1$ in the rest of the proof.

Fix a Riemannian metric on $Z$. We will denote by $\mu$ the corresponding
Riemannian measure on $Z$. 
Because of the transversality assumption and compactness of the divisor
$\D$, we may find positive constants
$c$ and $C$ such that 
$$c < |df(\xi)| < C\tag4-9$$
holds for all $\xi$ in some open neighborhood $W$ of $\D$.
Consider a smooth vector field $X$ on $Z$ which is orthogonal to $\D$ on $\D$
and such that the vectors of the field $X(\xi)$ have unit 
length for $\xi\in W$. Since $\<df,X\>=\pm |df|$ on $\D$, we may assume that the
neighborhood $W$ is so small that 
$$c/2 < |\<df(\xi),X(\xi)\>| <C\tag4-10$$ 
holds for all points in $W$. 

The neighborhood $W$ is fibred by integral trajectories of the field $X$
starting from the points of $\D$. If a point $x\in W$ belongs to a trajectory 
$x(t)$ starting at a point $p\in \D$ then using 
$\frac{d}{dt} f(x(t)) = \<df,X\>$
and (4-10) we obtain 
$$c\rho(x)/2 \le |f(x)| \le C\rho(x),\tag4-11$$
where $\rho(x)$ denotes the Riemannian distance between $x$ and $p$ along the
trajectory.

Thus we obtain that for small $\lambda$ the set $U_\lambda$ contains 
$\{\xi\in W; \rho(\xi) \le \lambda C^{-1}\}$; similarly, for small
$\lambda$ the set $U_\lambda$ 
is contained
in $\{\xi\in W; \rho(\xi) \le  2\lambda c^{-1}\}$. Therefore 
$$\lambda AC^{-1} \le F_\mu(\lambda) \le 2\lambda A c^{-1},\tag4-12$$
where $A$ denotes the area of $\D$. 

(4-12) shows that 
$F_\mu(\lambda) \sim \lambda$. \qed

\enddemo

When the divisor is not a submanifold the capacity may be not integral and
the spectral density function may be not equivalent to a power 
$\lambda^{\alpha}$. We consider a few illustrating examples.

\subheading{4.7. Example: Transversal crossing} Consider the 2-dimensional 
torus 
$T= \R/\Z\times \R/\Z$ with coordinates $x$ and $y$, defined modulo 
integers. Let $f: T\to T$ be the function
$f(x,y) = xy$. Consider the following torsion object 
$\X = (M_f:L^2(T)\to L^2(T))$. Then the divisor $\D(\X)$, which is equal
to $\{(x,y); x=0 \quad \text{or}\quad y=0\}$, has a singular
point $x=0, y=0$ (the cross). An easy computation of the spectral
density function of $\X$ with respect to the Lebesgue measure $\mu$
gives $F_\mu(\lambda) \sim \lambda(1-\log \lambda)$ which has capacity 1,
but {\it it is not dilatationally equivalent to} $\sim \lambda$.

\subheading{4.8. Example: Tangency} Let $Z = [-1,1]\times [-1,1]$ be the 
two-dimensional square with the coordinates $x$ and $y$ supplied with 
the Lebesgue measure. Let $k$ be a non-negative integer. 
Consider the function $f:Z\to \R$ given by $f(x,y) = y(y-x^k)$. 
As above, it defines the following torsion object 
$\X = (M_f:L^2(Z)\to L^2(Z))$. The divisor $\D(\X)$ is the union of
the interval $[-1,1]$ of the real axis and the parabola $y=x^k$; the 
divisor has a singular point at the origin. 
A computation of the capacity $\c(\X)$ in this example gives 
$$\c(\X) = \frac{2k}{k+1}.\tag4-13$$
This shows that the order of tangency $k$ can be recovered from the capacity
$\c(\X)$. Note also that in this example $1<\c(\X)<2$ for $k>1$ and 
$\c(\X)\to 2$ when $k\to\infty$.

\subheading{4.9. Example: Divisors of higher codimension} If the divisor
$\D(\X)$ is a submanifold but of codimension greater than 1, then the 
capacity $\c(\X)$ may be rational and not integral. 

Consider, for example the $n$-dimensional
torus $T^n$ with coordinates $x_1, x_2, \dots, x_n$ considered modulo $\Z$
and the function $f(x) = (\sum_{i=1}^nx_i^2)^m$, where $m>0$ is an integer. 
Then the divisor of the corresponding torsion object $\X_f$ consists of
one point $0$, and the Novikov - Shubin capacity $\c(\X_f)$ equals
$2m/n$. Moreover, the spectral density function is the 
power $\sim \lambda^{n/2m}$.

Now we formulate the following generalization of Proposition 4.3.

\proclaim{4.10. Theorem} Let $Z$ be a compact Hausdorff 
space supplied with a positive Radon measure $\mu$ such that there are no
nonempty open sets $U\subset Z$ with $\mu(U)=0$. 
Let $\E^i$, where $i=0, 1,2,\dots, N$, be a finite sequence of 
Hermitian vector bundles over $Z$. Suppose that for each $i$ there is given
a continuous bundle map $\partial^i:\E^i\to\E^{i+1}$ such that the sequence
$$(\E, \partial) = 
(0\to\E^0@>{\partial^0}>>\E^1@>{\partial^1}>>\E^2\dots @>{\partial^{(N-1)}}>>\E^N\to 0)\tag4-14$$
represents a cochain complex, i.e. $\partial^{i+1}\circ \partial^i = 0$. Denote by
${\Cal N}(\E,\partial)$ the set of all points $\xi\in Z$, such that the complex
$(\E_\xi,\partial_\xi)$ of the fibers over $\xi$ is not acyclic. Suppose that 
${\Cal N}(\E,\partial)$ has $\mu$-measure
zero. Then the extended cohomology $\H^\ast(L^2(\E),\partial^\ast)$ 
of the induced complex of $L^2$-sections
$$(L^2(\E),\partial^\ast) = 
(0\to L^2(\E^0)@>{{\partial^0}^\ast}>> L^2(\E^1)@>{{\partial^1}^\ast}>>\dots 
@>{{\partial^{(N-1)}}^\ast}>>L^2(\E^N)\to 0)\tag4-15$$
(understood as a graded object of the extended category $\eca$)     
is torsion in all dimensions and the union of the divisors
$$\bigcup_{i=0}^N \D(\H^i(L^2(\E),\partial^\ast))\tag4-16$$ 
coincides with
${\Cal N}(\E,\partial)$. Moreover, the extended cohomology 
$\H^\ast(L^2(\E),\partial^\ast)$ vanishes if and only if the set 
${\Cal N}(\E,\partial)$ is empty.\endproclaim 

\demo{Proof} We will deduce this statement from Proposition 4.3.

First we observe that the set ${\Cal N}(\E,\partial)$ is obviously
closed.
If $\xi\notin {\Cal N}(\E,\partial)$ then we can find a compact neighborhood $U$
of $\xi$ such that $U\cap {\Cal N}(\E,\partial) = \emptyset$ and then the complex
of bundles $(\E|_U,\partial|_U)$ is acyclic. Moreover, using induction we can
construct step by step a null-homotopy $S^i:\E^i|_U\to \E^{i-1}|_U$ by 
means of continuous bundle maps $S^i$ such that $S\partial+\partial S=\id$. 
Therefore, the corresponding complex of $L^2$-sections, constructed out of 
the complex $(\E|_U,\partial|_U)$,
is acyclic if considered as complex in the extended category. Thus, we 
obtain that for any measure $\mu'\ll \mu$ with $\supp(\mu')\subset U$
the spectral density functions of extended cohomology of complex (4-15)
vanish. Since the complement $Z - {\Cal N}(\E,\partial)$ is of full measure,
this shows that the extended cohomology of (4-15) is torsion and the
union of the divisors (4-16) is contained in ${\Cal N}(\E,\partial)$.

Suppose now that $\xi\in {\Cal N}(\E,\partial)$. Then the complex of fibers
$(\E_\xi,\partial_\xi)$ over $\xi$ has non-trivial cohomology in some dimension $i$. 
Consider the $i$-dimensional Laplacian: $\Delta^i:\E^i\to\E^i$. Here
$\Delta^i$ denotes $\hat \partial^i\circ \partial^i + \partial^{i-1}\circ \hat \partial^{i-1}$, where 
the sign {\it "hat"} denotes the adjoint bundle map with respect to the 
given Hermitian structures. $\Delta^i$ is a continuous bundle map 
$\E^i\to \E^i$,
and so we may apply Proposition 4.3 to it. Since at point $\xi$ the complex
of the fibers $(\E_\xi,\partial_\xi)$ has nontrivial cohomology in dimension $i$,
it follows that the $i$-dimensional Laplacian at the point $\xi$ 
(which equals to evaluation of $\Delta^i$ at $\xi$) has a non-trivial kernel. 
By Proposition
4.3 we obtain that then for any neighborhood $U$ of $\xi$ one may find a
measure $\mu'\ll\mu$ with $\supp(\mu')\subset U$ such that the spectral
density function of the torsion object 
$$({\Delta^i}^\ast: L^2(\E^i) \to L^2(\E^i))\tag4-17$$
is strictly positive for positive $\lambda$. But then formula (3-30) implies
that $\xi$ belongs to the divisor of the torsion object
$\T(\H^j(L^2(\E),\partial^\ast))$ for $j$ equal $i$ or $i+1$. \qed
\enddemo

In the analytic situation we can strengthen Theorem 4.10 by abandoning 
the assumption that the complex of bundles $(\E,\partial)$ is acyclic almost 
everywhere:

\proclaim{4.11. Theorem} Let $Z$ be a compact real analytic manifold, 
and let 
$$(\E, \partial) = 
(0\to\E^0@>{\partial^0}>>\E^1@>{\partial^1}>>\E^2\dots @>{\partial^{(N-1)}}>>\E^N\to 0)\tag4-18$$
be a cochain complex of real analytic Hermitian vector bundles over $Z$ 
and real analytic
bundle maps between them. For any point $\xi\in Z$ consider the complex 
$(\E_\xi,\partial_\xi)$ of the fibers above $\xi$ and denote by 
$\beta^i(\xi)$ the $i$-dimensional Betti number of this chain complex.
Then there is a minimal subset ${\Cal N}(\E,\partial)\subset Z$ 
of measure zero
such that all the functions $\beta^i(\xi)$, $i=0,1,\dots, N$ are constant
for $\xi\in Z - {\Cal N}(\E,\partial)$. We will denote by $\beta^i$ the
common value of $\beta^i(\xi)$ for $\xi\in Z - {\Cal N}(\E,\partial)$;
we will refer to it as to generic Betti number.

Fix a smooth measure $\mu$ on $Z$ coming from a Riemannian metric, and
consider the complex of $L^2$ sections $(L^2(\E),\partial^\ast)$:
$$
(0\to L^2(\E^0,\mu)@>{{\partial^0}^\ast}>> L^2(\E^1,\mu)
@>{{\partial^1}^\ast}>>\dots 
@>{{\partial^{(N-1)}}^\ast}>>L^2(\E^N,\mu)\to 0)
\tag4-19$$
which can be viewed as a projective chain complex in $\ca$.
Then the extended cohomology 
$\H^\ast(L^2(\E),\partial^\ast)$ of this complex has the following properties:
\roster
\item the von Neumann dimension of the projective part of the extended
cohomology of complex (4-19) is given by
$$\dim_{\tr} P(\H^i(L^2(\E),\partial^\ast)) = \mu(Z)\beta^i,\tag4-20$$ 
where $\tr$ denotes the trace on the von Neumann category $\ca$ constructed
by means of measure $\mu$ (cf. 2.7);
\item the union of the divisors of the torsion part of the extended
cohomology is given by 
$$\bigcup_{i=0}^N \D(\T(\H^i(L^2(\E),\partial^\ast))) = 
{\Cal N}(\E,\partial);\tag4-21$$
\item the extended cohomology of the $L^2$ complex (4-19) vanishes if and
only if the complex $(\E_\xi,\partial_\xi)$ is acyclic for all $\xi\in Z$.
\endroster
\endproclaim
\demo{Proof} Consider first a special case of Theorem 4.11, when it can
be trivially deduced from Theorem 4.10. 

Suppose that the given chain complex $(\E,\partial)$ over $Z$ can be
represented as a direct sum $(\E_0,\partial)\oplus(\E_+,\partial)$, where 
$(\E_0,\partial)$ and $(\E_+,\partial)$ are two chain complexes over $Z$ 
formed
by vector bundles and bundle maps, such that the first complex 
$(\E_0,\partial)$ is acyclic almost everywhere on $Z$ (as in Theorem 4.10)
and the second complex of vector bundles $(\E_+,\partial)$ has trivial 
differential. In this case we will say that the original chain complex 
$(\E,\partial)$ {\it splits}. 

It is clear that in the split case Theorem 4.11 is true. 
In fact, for any $i$ the rank of $\E_+^i$ equals to {\it the generic Betti 
number}
$\beta^i$; also the
complex of $L^2$ sections (4-19) is the direct sum 
$(L^2(\E_0),\partial^\ast)\oplus (L^2(\E_+),\partial^\ast)$ and so
the extended cohomology of $(L^2(\E),\partial^\ast)$ is the direct sum of the 
extended cohomology of $(L^2(\E_0),\partial^\ast)$ and 
$(L^2(\E_+),\partial^\ast)$. We may apply Theorem 4.10 to compute the 
extended cohomology of $(L^2(\E_0),\partial^\ast)$. The extended cohomology
of the second complex $(L^2(\E_+,\partial)$ in each dimension $j$ coincides 
with $L^2(\E_+^j)$. 

Let us show that one may reduce the general case to the split case by 
performing certain {\it blow up}. 

Our strategy will be as follows. 
Given a complex $(\E,\partial)$ of real analytic vector bundles as above, 
we will construct a compact Hausdorff space $\hat Z$ and a continuous map
$\phi:\hat Z \to Z$ such that
\roster
\item the induced complex of vector bundles $(\phi^\ast(\E),\partial)$ 
over $\hat Z$ splits (in the sense explained above);
\item $\phi$ is onto and for any point $p\in Z-{\Cal N}(\E,\partial)$ 
the preimage $\phi^{-1}(p)$ consists of precisely one point;
\item the set ${\Cal N}(\phi^\ast(\E),\partial)\subset \hat Z$ of 
{\it the special values} for the induced chain complex coincides with
$\phi^{-1}({\Cal N}(\E,\partial))$;
\item the interior of the set $\phi^{-1}({\Cal N}(\E,\partial))$ is
empty and there exists a positive Radon measure $\hat \mu$ on $\hat Z$ 
such that 
for any measurable subset $U\subset Z'$ we have $\mu'(U)=\mu(\phi(U))$.
\endroster
Since ${\Cal N}(\E,\partial)$ has measure zero (as a proper real analytic 
subset) we obtain that the canonical map
$L^2(\E^j,\mu)\to L^2(\phi^\ast(\E^j),\mu')$ is an isomorphism of Hilbert
$L^\infty(Z,\mu)$-modules and so we may identify the extended cohomology
of complex (4-19) with the extended cohomology of complex 
$(L^2(\phi^\ast(\E),\partial)$. The last complex splits, and so its extended
cohomology is given by the remark above.

We will now show how one can construct a blow up $\phi: \hat Z\to Z$ with the
above properties.

For each $i$ let $G_{\beta^i}(\E^i)\to Z$ denote the bundle of 
$\beta^i$-dimensional
Grassmannians associated to the bundle $\E^i$. 
Let $\beta$ denote the vector of generic Betti 
numbers 
$\beta = (\beta^0,\beta^1,\dots,\beta^N)$ and let $G_{\beta}(\E)\to Z$ denote the
fiber product
$$G_{\beta}(\E) = G_{\beta^0}(\E^0)\times_Z G_{\beta^1}(\E^1)\times_Z
\times\dots \times_Z G_{\beta^N}(\E^N).\tag4-22$$
A point $L\in G_{\beta}(\E)$ above some point $\xi\in Z$ is represented 
by a sequence $L=(L^0,L^1,\dots,L^N)$, where $L^i$ is a $\beta^i$-dimensional
subspace in the fiber $\E^i_\xi$ above $\xi$. We will denote by $\pi$ the
projection $G_{\beta}(\E)\to Z$. Thus, $\xi=\pi(L)$ in the notations above.

Denote by $Z'$ the subset of $G_{\beta}(\E)$ consisting of all points
$L=(L^0,L^1,\dots,L^N)\in G_\beta(\E)$ such that:
\roster
\item"{(a)}" $\partial_\xi(L^i)=0$ for all $i=0,1,\dots, N-1,$ where
$\xi=\pi(L)$;
\item"{(b)}" $\delta_\xi(L^i)=0$ for all $i=1,2,\dots, N.$ 
\endroster
Here $\delta_\xi$ denotes the adjoint of $\partial_\xi$ with respect to the
given Hermitian metrics on the fibers above $\xi$. The condition (a) 
means that the subspace $L^i$ is contained in $\ker(\partial_\xi)$ and
similarly for (b).

Denote by $\phi:Z'\to Z$ the restriction of $\pi$.

We see that $Z'$ is compact as a closed subset of the Grassmannian bundle 
$G_\beta(\E)$. 

If $\xi$ is a generic point of $Z$ (i.e. $\xi\notin {\Cal N}(\E,\partial)$) 
then the preimage $\phi^{-1}(\xi)$ consists of precisely one point
$L=(L^0,L^1,\dots,L^N)$, where $L^i$ is the space of harmonic vectors in the 
fiber $\E^i_\xi$. 

Now we may finally define the blow up $\hat Z$ as the closure of 
$\phi^{-1}(Z-{\Cal N}(\E,\partial))$. We will also denote by $\phi$ the
restriction map $\phi|_{\hat Z}$.

The induced chain complex of bundles 
$(\phi^\ast(\E),\partial)$ over $\hat Z$ splits; its splitting over a point
$L\in \hat Z$ is given by assigning $(\phi^\ast(\E^i)_+)_L = L^i$ and
$(\phi^\ast(\E^i)_0)_L = {L^i}^\perp$. 

It is easy to see that the conditions (1) - (4) mentioned above are
satisfied.\qed

Note, that as follows from the construction, the bow up 
$\hat Z$ has structure of a real algebraic set. 
\enddemo

\subheading{4.12. Novikov - Shubin invariants and the  germ-cohomology} 
Assume that additionally to the hypotheses 
of Theorem 4.11, the manifold of parameters $Z$ is one-dimensional.
Then it is possible to compute explicitly
the Novikov-Shubin numbers of the extended cohomology localized at different
points of the divisor. We will see that
under these assumptions the capacity of the extended cohomology
is always an integer -- this should
be compared with examples 4.8
and 4.9, which show that the capacity may be not integral (but rational!) 
in the real analytic situation.

In fact, I believe that {\it in the real analytic situation of Theorem 4.11
the capacity of the extended cohomology of complex 
(4-15) is always rational.} Compare with conjecture 7.1 in \cite{LL}.

In the case of one-dimensional parameter the capacity
of the extended $L^2$ cohomology can be expressed using the notions of 
{\it germ-complex and germ-cohomology}, which were introduced in \cite{Fa2}
and \cite{FL}. These notions are defined in a more general 
situation when we have a family of elliptic 
complexes (or a finite dimensional family, like (4-14)) depending on a 
parameter $t\in (t_0-\epsilon,t_0+\epsilon)$ in a real analytic fashion. 
The germ complex
is constructed similarly to complex (4-15), but instead of $L^2$ sections,
one considers {\it the germs at $t_0$ of real analytic sections.} 
More precisely, denote by $\OO\E^i$ the set of germs of all real analytic
sections of (4-14), defined in a neibourhood of $t=t_0$. The differential
$\partial:\E^i\to \E^{i+1}$ naturally defined the map on the sections
$\partial:\OO\E^i\to\OO\E^{i+1}$.
This produces a complex $(\OO\E,\partial)$
(called the germ-complex) of free modules 
over the ring $\OO$ of germs at $t=t_0$ of real analytic 
complex valued functions of the parameter $t$. 
We refer the reader
to \cite{Fa2, FL} for more information and applications. 

The cohomology of the germ complex $H^i(\OO\E,\partial)$ is a finitely
generated module over $\OO$. Since $\OO$ is a discrete valuation ring, 
the module
$H^i(\OO\E,\partial)$ has well defined {\it $\OO$-torsion submodule}
$\tau^i\subset H^i(\OO\E,\partial)$; the whole module 
$H^i(\OO\E,\partial)$ is in fact a direct
sum of a free $\OO$-module (of rank equal to the generic Betti number 
$\beta^i$) and $\tau^i$. The torsion submodule is finite 
dimensional as a vector space over $\C$; it has a canonical
nilpotent endomorphism 
$\hat t:\tau^i\to\tau^i$ - multiplication by $t-t_0\in\OO$. The minimal  
number $n$ such that $\hat t^n$ annihilates $\tau^i$ (i.e. 
$\hat t^n\cdot \tau^i = 0$) will be called 
{\it the height of $\tau^i$}.

Now we may formulate an addition to Theorem 4.11 in the case
of one-dimensional parameter space.

\proclaim{4.13. Theorem} If manifold $Z$ in Theorem 4.11 is 
one-dimensional, then for any $i$ the divisor of the extended
cohomology $\H^i(L^2(\E),\partial^\ast)$ of complex (4-15) consists of 
finitely many isolated points.
If $t_0$ denotes one of the points of the divisor, then the capacity of 
the extended cohomology $\H^i(L^2(\E),\partial^\ast)$
with respect to a measure of the 
form $\nu=\chi_{[t_0-\epsilon', t_0+\epsilon']}(t)\times\mu$ (where 
$\chi$ denotes the characteristic function of the indicated
interval and $\mu$ is
the Lebesgue measure) for small enough $\epsilon'$ equals to
the height of the torsion part of the germ cohomology
$$\c_\nu(\H^i(L^2(\E),\partial^\ast)) = \text{height of}\quad \tau^i.$$
In particular, this capacity is integral.
\endproclaim
\demo{Proof} First of all, applying the recipe described in section 3.14, 
we realize that to compute the
spectral density function of $\H^i(L^2(\E),\partial^\ast)$ we have to 
consider the self-adjoint bundle map 
$\partial^\ast\partial:\E^{i-1}\to\E^{i-1}$. Using the Rellich - Kato theorem
(cf. \cite{K}, chapter 7), we may find {\it real analytic functions}
$$\lambda_1,\dots,\lambda_m: U\to\R,$$
where $m$ is the rank of the bundle $\E^{i-1}$ and $U$ is a neighbourhood
of $t_0$, such that for any $t\in U$ the numbers 
$\lambda_1(t),\dots,\lambda_m(t)$ are all the eigenvalues of 
$\partial^\ast\partial$. 
There can be several functions $\lambda_i(t)$, which are identically zero.
We will assume that the $\lambda_i(t)$'s are numerated such that 
$\lambda_i(t)\equiv 0$ for $1\le i<m'$ and for $m'\le i\le m$ the functions
$\lambda_i(t)$ are not identically zero.

Now, similarly to (4-4), the spectral density 
function
$F_\nu(\lambda)$ of $\H^i(L^2(\E),\partial^\ast)$ is given by
$$F_\nu(\lambda) = \sum_{i=m'}^m\nu\{t\in U; \lambda_i(t)\le \lambda^2\}.$$
Since $\lambda_i(t)$ is real analytic and non-negative around $t_0$, we
may write
$$\lambda_i(t) = (t-t_0)^{2k_i}\gamma_i(t),\qquad i= m',\dots, m,$$
where $k_i$ is a non-negative integer and $\gamma_i(t_0)\ne 0$. 
Note, that the point $t_0$ belongs to the divisor if and only if at least
one of these numbers $k_i$ is non-zero (by Theorem 4.11). Thus, we
obtain that the points of the divisor are isolated and also, since
$$\nu\{t\in U; \lambda_i(t)\le \lambda^2\}\sim \lambda^{k_i^{-1}},
\quad\text{for}\quad k_i\ne 0,$$ 
the Novikov - Shubin number of $\H^i(L^2(\E),\partial^\ast)$ is 
$\min\{k_i^{-1}\}$. Therefore, the capacity 
$\c_\nu(\H^i(L^2(\E),\partial^\ast))$ is $\max\{k_i\}$. 

The Rellich - Kato theorem \cite{K}, chapter 7, gives us also an orthonormal
system of eigensections $s_1(t), \dots, s_m(t)$ of $\E^{i-1}$ for
$t\in U$ such that $\partial^\ast\partial s_i(t)=\lambda_i(t)s_i(t)$ and
such that $s_i(t)$ are real analytic in $t$. For $m'\le i\le m$ 
the sections 
$\lambda_i^{-1/2}(t)\cdot \partial s_i(t)$ of $\E^i$
are then analytic and form an orthonormal basis
in the space of analytic sections of $\E^i$
which belong to the image of $\partial$
after multiplication by a power of $t-t_0$ (compare with the definition of 
the germ-complex in \cite{Fa2}). Thus, it follows that the height of 
torsion submodule $\tau^i$ of the cohomology of the germ-complex
also equals $\max\{k_i\}$.
\qed
\enddemo

\heading{\bf \S 5. Abelian extensions of categories of Hilbertian 
representations}\endheading 

In this section we mention briefly the similar construction of
building the abelian extension of an additive category in the case 
when instead of Hilbert representations we start with a category of 
{\it Hilbertian representations}. This situation appears to be more natural
from the point of view of many topological applications. For
example, Hilbertian representations are important if
one studies the $L^2$ torsion, combinatorial and analytic
(as in \cite{CFM}). Another class of examples (which seem to be important for
constructing a De Rham version of the extended $L^2$-cohomology, cf. \S 7) 
provide Sobolev spaces of sections of vector bundles with Hilbert fiber 
over a compact manifold. We will apply the results of this section in \S 7.

\subheading{5.1. Hilbertian spaces} Recall that {\it a Hilbertian space} 
(cf. \cite{P}) is a topological vector space $\H$ which
is isomorphic to a Hilbert space in the category of topological vector 
spaces. In other words, there exists a scalar product on $\H$ such that 
$\H$ with this scalar product is a Hilbert space with the originally given 
topology. Such scalar
products are called {\it admissible}. Given one admissible scalar product
$\< \ ,\ \>$ on $\H$, any other admissible scalar product is given by
$$\<x,y\>_1\ =\ \<Ax,y\>,\quad x,y\in \H,\tag5-1$$
where $A:\H\to \H$ is an invertible positive operator $A^\ast=A,\quad A>0$.

Hilbertian spaces naturally appear as Sobolev spaces of vector bundles,
cf. \cite{P}.

\subheading{5.2. Hilbertian representations} Let $\A$ be an algebra with an 
involution, as in 1.1. {\it A Hilbertian
representation of $\A$ } is a Hilbertian topological vector space $\H$
supplied with a left action of $\A$ by continuous linear maps 
$\A\to \L(\H,\H)$
such that there exists an admissible scalar product on $\H$ with respect to 
which the given action is a Hilbert representation, cf. 1.1. We will say
that this scalar product is {\it admissible with respect to the given
Hilbertian representation}. Again, given one such admissible scalar product, 
all other admissible scalar products can be constructed as in the previous 
paragraph, where $A$ is any invertible positive operator {\it commuting with 
the action of $\A$}.  

A morphism $f:\H_1\to \H_2$ between two Hilbertian representations is 
defined in the same way as in 1.1, i.e. as a bounded linear map commuting 
with the action of the algebra $\A$.

\subheading{5.3. Hilbertian categories} As in 1.1, we may consider a 
subcategory $\ca$ of the category of Hilbertian
representations of a given algebra with involution $\A$. 

We will say that
$\ca$ is a {\it Hilbertian category} if any object $\H$ of $\ca$ has a 
$\ca$-admissible scalar product (cf. below) and $\ca$ satisfies condition 
(i) of section 1.2 and also condition (ii) of section 1.2 in the following 
edition:

\roster
\item"{(ii$^\ast$)}" {\it For any $\H\in\ob(\ca)$ the dual representation
$\H^\ast$ also is an object of $\ca$ and for any morphism $\phi:\H_1\to\H_2$ 
of $\ca$ the adjoint operator $\phi^\ast:\H_2^\ast\to\H_1^\ast$ is also a 
morphism of $\ca$.}
\endroster

Note, that the dual representation $\H^\ast$ is defined as the space of all
bounded anti-linear functionals $\phi:\H\to\C$, 
$\phi(\lambda x)=\overline\lambda \phi(x)$ for $x\in\H$ and $\lambda\in \C$.
The $\A$-module structure on $\H^\ast$ is given by $(a\phi)(x)=\phi(a^\ast x)$
for $a\in\A$ and $x\in\H$. 

If $\ca$ is a category of Hilbertian representations 
and $\H\in\ob(\ca)$, then we can specify naturally a class of admissible 
scalar products on $\H$ by declaring an  
admissible scalar product $\<\ ,\ \>$ on $\H$ to be {\it $\ca$-admissible} 
if the isomorphism $\H\to \H^\ast$ determined by $\<\ ,\ \>$ belongs to $\ca$.

If $\<\ ,\ \>$ and $\<\ ,\ \>_1$ are two $\ca$-admissible scalar products on 
$\H$ then the corresponding operator $A$ (which appears in (5-1)) and its 
inverse $A^{-1}$ are morphisms of $\ca$.

\subheading{5.4. Abelian extension of a Hilbertian category}
Repeating construction described in section 1.3, for any Hilbertian category
$\ca$ one obtains an extended abelian category $\eca$ containing $\ca$ 
as a full subcategory. We will call $\eca$ {\it the abelian extension of} 
$\ca$. Any object of the extended category $\eca$ is
represented as a morphism $(\alpha:A'\to A)$ in $\ca$. As in section 1.8,
there is a full embedding $\ca \to\eca$; similarly to Proposition 1.9 we have
that any object of the original category $\ca$ is projective in $\eca$ 
and conversely, any projective object of $\eca$ is isomorphic to an object 
of $\ca$ inside $\eca$. 

It is easy to check that all the statements and the 
arguments of subsections 1.1 - 1.11 equally apply to the Hilbertian  
situation as well.

\subheading{5.5. Hilbertian von Neumann categories} A Hilbertian 
category $\ca$ will be called {\it von Neumann category} if it satisfies
condition (v) of section 2.1. 

An object $\H$ of a Hilbertian von Neumann category will be called {\it 
finite} if the only closed $\ca$-submodule $\H_1\subset\H$ which is 
isomorphic to $\H$ in $\ca$ is $\H_1=\H$.

A Hilbertian von Neumann category will be called {\it finite} iff all its  
objects are finite.

Examples of Hilbertian von Neumann categories can be obtained from examples
considered in \S 2 by forgetting the scalar product.

\subheading{5.6. Torsion objects of the extended category} Suppose that
$\ca$ is a finite Hilbertian von Neumann category. Then similarly to \S 3,
we may define the notion of a torsion object. Denoting by $\tca$ the full
subcategory of the extended category $\eca$ generated by torsion objects, 
we will have (similarly to Corollary 3.3) that the torsion subcategory 
$\tca$ is {\it an abelian subcategory} of the extended category.

\subheading{5.7} The main distinction between this Hilbertian
version and the situation of \S 2 is that in the Hilbertian case for any 
$\H\in\ob(\ca)$
we do not have a fixed $\ast$-operator on $\Hom_{\ca}(\H,\H)$, and so 
it is not a von Neumann algebra. Instead,
we have many $\ast$-operators (each corresponding to a choice of an 
admissible
scalar product on $\H$), and considered with any of these involutions the 
algebra $\Hom_{\ca}(\H,\H)$ becomes a von Neumann algebra. We observe that
the notion of positivity of an element of $\Hom_{\ca}(\H,\H)$ also 
depends on the choice of an admissible scalar product on $\H$.

\subheading{5.8. Traces} {\it A trace on a Hilbertian von Neumann category} 
$\ca$ is a function $\tr$ which assigns
to each object $\H\in\ob(\ca)$ a linear functional  
$$\tr_{\H}: \Hom_{\ca}(\H,\H)\to \C\tag5-2$$
such that for any pair of representations $\H_1, \H_2\in\ob(\ca)$ 
the corresponding 
traces $\tr_{\H_1}$ and $\tr_{\H_2}$ are compatible in the following sense: 
if $f\in\Hom_{\ca}(\H_1, \H_2)$ and $g\in \Hom_{\ca}(\H_2,\H_1)$
then 
$$\tr_{\H_1}(gf) = \tr_{\H_2}(fg).\tag5-3$$

We may call trace $\tr$ {\it non-negative} if $\tr_{\H}(e)$ is real and
non-negative for any
idempotent $e\in \Hom_{\ca}(\H,\H)$, $e^2=e$.

\subheading{5.9. von Neumann dimension and the Novikov-Shubin invariants}
Any trace on a von Neumann category defines {\it a von Neumann dimension 
function} (as in 2.8); the dimension is always non-negative and real valued 
iff the trace $\tr$ is non-negative (in the above sense).

Using a non-negative trace on $\ca$ one defines for any torsion 
object $\X$ of 
the extended category $\eca$ {\it the spectral density function} (as in 3.7) 
and then {\it the Novikov-Shubin number} $\ns(\X)$ or the {\it capacity}
$\c(\X)$ (as in 3.9). 

Hilbertian versions of Proposition 3.10 and Corollary 3.11 hold.

The following Lemma justifies our definition of the trace
on Hilbertian category.

\proclaim{5.10. Lemma} Suppose that $\ca$ is a Hilbertian category and
$\H\in\ob(\ca)$. Let $\tr_{\H}: \Hom_{\ca}(\H,\H) \to \C$ be a linear
functional which is traceful, i.e. $\tr_{\H}(fg)= \tr_{\H}(gf)$ for any
$f,g\in \Hom_{\ca}(\H,\H)$. Suppose that $\tr_H$ is non-negative and 
normal with respect to one choice of a $\ca$-admissible scalar product on 
$\H$ (i.e. it assumes real non-negative values on the positive cone
$\Hom_{\ca}(\H,\H)^+$ taken with respect to the given $\ast$-operator). 
Then it is
non-negative and normal with respect to any other choice of a 
$\ca$-admissible scalar product on $\H$ .\endproclaim

\demo{Proof} This follows from \cite{Di}, part I, chapter 6, Proposition 1
using the fact that any idempotent in a von Neumann algebra is conjugate
to a projection.\qed
\enddemo

\heading{\bf \S 6. Extended $L^2$ cohomology of cell complexes}\endheading

Our aim in this section is to construct homology and cohomology theories on 
the category of finite polyhedra with values in the extended abelian 
category $\eca$. Then we compute the extended cohomology of mapping tori.

We will start with some very general remarks, which will be
then used in the main construction, cf. 6.5. Our goal is to describe {\it
the coefficient systems}.

\subheading{6.1} Let ${\Cal C}$ be an additive category. 
The objects of ${\Cal C}$ will be denoted $\X$, $\Y$ etc., and the
set of morphisms between $\X$ and $\Y$ will be denoted 
by $\hom_{\Cal C}(\X,\Y)$.

Let $\Lambda$ be a ring with a unit. We will assume that $\Lambda$ has 
IBN (invariance of the basis number) property, cf. \cite{C}. For our 
applications we need only the 
case $\Lambda = \C[\pi]$, where $\pi$ is a discrete group; then this property
is automatically satisfied, cf. \cite{C}.

\subheading{Definition} {\it A structure of left $\Lambda$-module} on an 
object $\X\in\ob({\Cal C})$ consists
in specifying a ring homomorphism $\rho:\Lambda\to \hom_{\Cal C}(\X,\X)$.
Similarly, a structure of right $\Lambda$-module on $\X$ is given by a 
ring homomorphism  $\rho:\Lambda^{op}\to \hom_{\Cal C}(\X,\X)$, where
$\Lambda^{op}$ is the opposite ring of $\Lambda$.

If one of these structures is chosen, we will say that {\it $\X$ is a left 
(right) $\Lambda$-module} in ${\Cal C}$. We will usually omit $\rho$ from 
the notation. The left (right) $\Lambda$-modules in the category of abelian 
groups correspond to the usual notion of left (right) $\Lambda$-module.

If $\X$ and $\Y$ are two left (right) $\Lambda$-modules in ${\Cal C}$ then 
a morphism $f\in\hom_{\Cal C}(\X,\Y)$ is called a 
{\it $\Lambda$-homomorphism} if 
for any $a\in \Lambda$ the following diagram commutes
$$
\CD
\X@>f>>\Y\\
@V{\rho(a)}VV   @VV{\rho(a)}V\\
\X@>f>>\Y
\endCD\tag6-1
$$
Thus we have a category of left $\Lambda$-modules in a given additive 
category ${\Cal C}$, which we will denote $\Lambda$-mod-${\Cal C}$; the 
similar category of right $\Lambda$-modules in ${\Cal C}$ will be denotes
${\Cal C}$-mod-$\Lambda$.

If $\X$ is a left $\Lambda$-module in ${\Cal C}$ then for any 
$\Y\in\ob({\Cal C})$
the set of morphisms $\hom_{\Cal C}(\Y,\X)$ has a natural structure of a 
left $\Lambda$-module, where for $a\in \Lambda$ and 
$f\in \hom_{\Cal C}(\Y,\X)$ we define
$$a\cdot f = \rho(a) \circ f : \Y \to \X.\tag6-2$$
Thus, if $\X$ is a $\Lambda$-module in ${\Cal C}$, 
the functor $\hom_{\Cal C}(\cdot,\X)$ assumes its values in the 
category of left $\Lambda$-modules.

Let $F$ be a free finitely generated left $\Lambda$-module (in the usual 
sense) and let $\X$ be a left $\Lambda$-module in an additive category 
${\Cal C}$. We are now going to define 
$$\om_{\Lambda}(F,\X),\tag6-3$$
which will be an object of category ${\Cal C}$. It will depend functorially
on both $F$ and $\X$, as stated in the following Proposition.

\proclaim{6.2. Proposition} There exists a bifunctor 
$$\om_{\Lambda}(\cdot,\cdot): \ \{\text{free f.g. left
$\Lambda$-modules}\}\times (\Lambda\text{-mod-}{\Cal C}) \ 
\to \ {\Cal C},\tag6-4$$
which is contravariant with respect to the first variable and covariant with
respect to the second variable and such that there is a natural isomorphism
$$\hom_{\Cal C}(\Y,\om_\Lambda(F,\X)) \simeq  
\Hom_\Lambda(F, \hom_{\Cal C}(\Y,\X))\tag6-5$$
where $\Y\in\ob({\Cal C})$, $\X\in\ob(\Lambda\text{-mod-}{\Cal C})$ and  
$F$ is a free finitely generated left $\Lambda$-module. If the category
${\Cal C}$ is abelian, then (6-4) is exact as a functor of $\X$.
\endproclaim
\demo{Proof} For each free f.g. left $\Lambda$-module $F$ choose a base
$e = (e_1, e_2, \dots, e_n)$, where $n=\rank F$. Here we should assume that
we are dealing with a small category of free finitely generated left
$\Lambda$-modules and then such choice is possible. The functor which
we will construct, considered up to isomorphism, is independent on 
this choice.  

Define 
$\om_{\Lambda}(F,\X)$ as $\X\oplus\X\oplus\dots\oplus\X$ ($n$ times).
It obviously behaves functorially with respect to $\X$ and is exact if 
${\Cal C}$ is abelian.
If $F'$ is another free f.g. left $\Lambda$-module with the base
$e' = (e'_1, e'_2, \dots, e'_m)$, where $m=\rank F'$, then any 
$\Lambda$-homomorphism $\phi:F\to F'$ is represented in the chosen bases
by a matrix $(a_{ij})$, where $\phi(e_i)=\sum_{j=1}^ma_{ij}e'_j$, with
$a_{ij}\in\Lambda$. Then we define $\phi^\ast: \om_{\Lambda}(F',\X)\to
\om_{\Lambda}(F,\X)$ as the morphism which maps the $j$-th copy of $\X$
to $i$-th copy of $\X$ by $\rho(a_{ij})$. \qed
\enddemo

\subheading{6.3. Change of rings} Suppose that $\phi:\Lambda\to\Lambda'$ is
a ring homomorphism. Then any left (right) $\Lambda'$-module $\X$ in 
${\Cal C}$ determines via $\phi$ a left (right) $\Lambda$-module in 
${\Cal C}$ which we denote $\phi^\ast(\X)$.

Also, suppose that $F$ and $F'$ are free finitely generated modules over
$\Lambda$ and $\Lambda'$ correspondingly and let $f:F\to F'$ be a
$\Lambda$-homomorphism (via $\phi$). Then $f$ induces canonically the morphism
$$f^\ast:\om_{\Lambda'}(F',\X) \to \om_{\Lambda}(F,\phi^\ast(\X)).\tag6-6$$ 
It is defined by using the definition of $\om_{\Lambda'}(F,\X)$
above:  free basis of $F$ and $F'$ give a representation of the 
map $f$ by a matrix with entries in $\Lambda'$ which then may act
(using the module structure) on the direct sums of copies of $\X$
similarly to the arguments above.

We will mention here one situation when homomorphism (6-6) is an isomorphism;
we will refer to it later. Suppose that $\X$ is a left $\Lambda$-module 
in ${\Cal C}$ and $\rho:\Lambda\to
\hom_{\Cal C}(\X,\X)$ is the action homomorphism. Its kernel $\ker(\rho)$
is an ideal in $\Lambda$. Let $\Lambda'$ be the factor-ring 
$\Lambda/\ker(\rho)$. $\X$ is well defined as a $\Lambda'$-module. For any free
finitely generated left $\Lambda$-module $F$ consider the morphism
$f:F\to  F'=F/\ker(\rho)F$. Then the corresponding morphism (6-6) is an isomorphism.

\subheading{6.4. Tensor products} We will consider here a similar construction
in the covariant modification. 

Suppose that $F$ is a {\it left} 
free f.g. $\Lambda$-module and $\X$ is a {\it right} $\Lambda$-module in an
additive category ${\Cal C}$. Similarly to Proposition 6.2, one defines
the tensor product $\X\ts F$, giving a bifunctor 
$$\ts : \ {\Cal C}\text{-mod-}\Lambda \times 
\{\text{free f.g. left $\Lambda$-modules}\}
\to \ {\Cal C},\tag6-7$$
which is covariant with respect to both variables. For any choice of a free
basis of $F$, the tensor product $\X\ts F$ can be identified with the direct
sum $\X\oplus\dots\oplus \X$ ($n$ copies), where $n$ is the rank of $F$.

If $\Y\in\ob({\Cal C})$, then the set $\hom_{\Cal C}(\Y,\X)$ has a structure
of a right $\Lambda$-module (in the usual sense) and there is 
natural isomorphism
$$\hom_{\Cal C}(\Y,\X\ts F) \simeq 
\hom_{\Cal C}(\Y,\X)\otimes_\Lambda F,\tag6-8$$
similar to (6-5).

\subheading{6.5. The basic construction} Now we may define extended homology
and cohomology of cell complexes. The construction in the form presented here,
generalizes the construction of \cite{Fa1}, \S 6. 

Let $K$ be a connected CW complex with a base point. We will assume that 
$K$ has finitely 
many cells in each dimension. Denote by $\pi=\pi_1(K)$ the fundamental group. 
Let $\tilde K$ denotes the universal cover
of $K$. We will suppose that a base point has been chosen in $\tilde K$,  lying
above the base point of $K$. Then the group of covering transformations
of $\tilde K$ can be identified with $\pi$.
Consider the chain complex $C_\ast(\tilde K)$ constructed using the cells
of the universal covering $\tilde K$, which are the lifts of the cells of $K$. 
The fundamental group
$\pi$ acts on  $C_\ast(\tilde K)$ from the left and in each dimension $i$
the chain module 
$C_i(\tilde K)$ is a free $\C[\pi]$-module with the number of generators
equal to the number of $i$-dimensional cells in $K$.

Let $\A$ be a $\ast$-algebra and let $\ca$ be a Hilbert category of 
$\ast$-representations of $\A$, cf. 1.2. (Later we will impose an
additional assumption that $\ca$ is
a finite von Neumann category, but most of the constructions hold in this 
more general context.) Let $\eca$ denote the extended 
abelian category of $\ca$, cf. \S 1. 

Assume that we are given a left $\C[\pi]$-module $\M$ in
$\eca$, as defined in section 6.1. We will consider some examples later in 
this section. The module $\M$ will play role of {\it a coefficient system}. 
Let us emphasize that in our present approach $\M$ may equally be projective, 
torsion, or a combination of a torsion module and  a projective. This 
generality may look superfluous, but (as we are going to show in another 
place) only having the opportunity to use torsion Hilbert spaces as
coefficient systems one may generalize naturally and fully many classical 
techniques, for example, the spectral sequences.

Applying the functor $\om_\Lambda(\cdot, \M)$, where 
$\Lambda=\C[\pi]$, to the chain complex
$C_i(\tilde K)$, we obtain the following cochain complex 
$\om_{\Lambda}(C_\ast(\tilde K),\M)$ in $\eca$:
$$\dots\leftarrow \om_\Lambda(C_{i+1}(\tilde K), \M)  
\leftarrow \om_\Lambda(C_{i}(\tilde K), \M) 
\leftarrow \om_\Lambda(C_{i-1}(\tilde K), \M)\leftarrow\dots\tag6-9$$  
The cohomology of this complex, understood as an object of the extended
category $\eca$, will be denoted 
$$\H^i(K,\M) = H^i(\om_{\Lambda}(C_\ast(\tilde K),\M));\tag6-10$$
it will be 
called {\it  extended cohomology of $K$ with coefficients in $\M$}. 

Intuitively, we may view the  cochains $c\in \om_\Lambda(C_{\ast}(\tilde K), \M)$ 
as functions which assign to cells of $\tilde K$  {\it "elements of $\M$ "} such
that 
$$c(ge)=\rho(g)(c(e))\tag6-11$$ 
for all cells $e$ of $\tilde K$ and all $g\in\pi$. Here 
$\rho:\Lambda\to \End_{\eca}(\M)$ is the module action.

Now we will define the similar notion of extended homology.  
Here we will assume that we are given a {\it right} $\C[\pi]$-module $\M$
in the extended category $\eca$.
We define
$$\H_i(K,\M) = H_i(\M\ts C_\ast(\tilde K)),\tag6-12$$
{\it the extended homology of $K$ with coefficients in $\M$}.
Here $\ts$ denotes the tensor product introduced in section 6.4.

Note that the extended homology and  cohomology are objects of the 
extended category and may have non-trivial torsion parts
even if the module $\M$ was projective.

\subheading{6.6} We will consider here some simple properties of the extended
cohomology considered as a functor of $K$. Similar properties are valid for
the extended homology.

Let $f:K\to K'$ be a continuous map preserving the base points between CW 
complexes having finitely many cells in each dimension. 
It induces a homomorphism
$f_\ast:\pi \to \pi'$, where $\pi = \pi(K)$ and $\pi' = \pi(K')$ are 
the fundamental groups. If $\M$ is a 
$\C[\pi']$-module in $\eca$, then the homomorphism $f_\ast$ allows 
to view $\M$ as a $\C[\pi]$-module in category $\eca$, which we will 
denote $f^\ast\M$. Now, the map $f$ lifts uniquely to a $\pi$-equivariant 
map $\tilde f:\tilde K\to\tilde K'$ preserving the base points. 
$\tilde f$ induces
a chain map $C_\ast(\tilde f): C_\ast(\tilde K)\to C_\ast(\tilde K')$. 
Applying to 
the last map the change of rings morphism (6-6), we obtain a chain map 
$$\om_{\Lambda'}(C_\ast(\tilde K'), \M) \to 
\om_\Lambda(C_\ast(\tilde K), f^\ast\M);\tag6-13$$
it induces the morphism of the extended cohomology
$$f^\ast: \H^\ast(K',\M) \to \H^\ast(K,f^\ast\M).\tag6-14$$
The induced morphism (6-14) depends only on the homotopy class of $f$, 
since the chain homotopy class of the map 
$C_\ast(\tilde f)$ (as it is well known) depends only on the homotopy class
of $f$. 

\proclaim{6.7. Corollary} The extended cohomology $\H^\ast(K,\M)$ is a
homotopy invariant of the CW complex $K$.\endproclaim

If the category $\ca$ is a finite von Neumann category, then to any trace
on $\ca$ we associate the Novikov - Shubin numbers (cf. 3.9) in order to
measure the size of the torsion part of the 
extended homology and cohomology. This method
gives {\it homotopy invariants o}f $K$ by Corollary 6.7 above - this fact
was proven first by M.Gromov and M.Shubin \cite{GS}.

\proclaim{6.8. Mayer - Vietoris sequence} Suppose that $K$ is a connected
CW complex with base point having finitely many cells in each dimension
and $\M$ is a $\C[\pi]$-module in $\eca$, where $\pi=\pi(K)$. Let $K_1$ and
$K_2$ be two connected subcomplexes containing the base point such that
$K = K_1 \cup K_2$ and the CW complex $K_0 = K_1\cap K_2$ is connected.
Denote $\M_i=\M|_{K_i}$, where
$i=0,1,2$. Then the following sequence in category $\eca$ is exact
$$\dots \to \H^i(K,\M)\to \H^i(K_1,\M_1)\oplus\H^i(K_2,\M_2)\to
\H^i(K_0,\M_0)@>\delta>>\H^{i+1}(K,\M)\dots\tag6-15$$
The morphisms of this exact sequence (except $\delta$) are induced by the
inclusions as usual.
\endproclaim
The proof repeats the standard arguments and will be skipped.

Now we will summarize some properties of the extended cohomology as function 
of $\M$.

\proclaim{6.9. Proposition} The extended cohomology $\H^i(K,\M)$ is a
covariant functor of $\M$. For any exact sequence of $\C[\pi]$-modules in
category $\eca$ 
$$0 \to \M' \to \M \to \M''\to 0\tag6-16$$
there is a natural long exact sequence
$$\dots \H^i(K,\M')\to\H^i(K,\M)\to\H^i(K,\M'')\to 
\H^{i+1}(K,\M')\to\dots\tag6-17$$
in $\eca$.
If $\ca$ is a finite von Neumann category (cf. \S 2) and $\M$ is torsion 
(as an object
of $\eca$) then all extended cohomology $\H^\ast(K,\M)$ is torsion.
\endproclaim
\demo{Proof} Since the functor $\om_{\Lambda}(F,\M)$ is exact with respect
to $\M$ (cf. Proposition 6.2), we obtain the following exact sequence
of cochain complexes
$$0\to \om_{\Lambda}(C_\ast(\tilde K),\M')\to 
\om_{\Lambda}(C_\ast(\tilde K),\M)\to \om_{\Lambda}(C_\ast(\tilde K),\M'')
\to 0,\tag6-18$$
which produces by the general laws of homological algebra
the exact sequence (6-17). If $\M$ is a torsion object then
all objects $\om_{\Lambda}(C_i(\tilde K),\M)$ are torsion and hence the
cohomology is torsion (by corollary 3.3). \qed\enddemo

Now we are going to discuss some examples of $\C[\pi]$-modules in the 
abelian extensions of Hilbert and von Neumann categories,
and the corresponding extended cohomology of polyhedra.

\subheading{6.10. Example: The regular representation}  
The simplest example of a $\C[\pi]$-module in $\eca$ can be obtained 
if we are given a representation of $\pi$ acting on
a Hilbert space $\H$, such that for any $g\in\pi$ the morphism $\H\to\H$
of multiplication by $g$, belongs to a Hilbert category $\ca$. Then $\H$
is a $\C[\pi]$-module in $\eca$.

The most familiar example of this sort is given by the regular representation
$\ell^2(\pi)$ with $\pi$ acting from the left by translations. 
The corresponding category $\ca$ in this case is the finite
von Neumann category
described in example 5 of section 2.6. The algebra
$\A$  is ${\Cal N}(\pi)^\bullet$, i.e. the algebra opposite to 
${\Cal N}(\pi)$, the von Neumann algebra of $\pi$, cf. 2.6, example 5.

The cochains in this situation have very transparent geometric meaning.
Namely, any cochain $c\in\om_{\Lambda}(C_\ast(\tilde K),\ell^2(\pi))$ is a
function which assigns an element $c(e)\in \ell^2(\pi)$ to any cell $e$
of the universal covering $\tilde K$ such that $c(ge)=gc(e)$ holds 
(cf. (6-11)) for
any $g\in \pi$. Writing $c(e)=\sum_{g\in\pi}\alpha_g\cdot g$, 
where $\alpha_g\in\C$ satisfy $\sum |\alpha_g|^2 < \infty$, we may
consider a cochain on $\tilde K$ with values in $\C$, where
$h(e)=\alpha_1$ (here $1$ denotes the unit element of $\pi$). 
Then $h(g^{-1}e)=\alpha_g$ and we obtain that $h$ is a
$\C$-valued cochain {\it satisfying the $L^2$-condition}:
$$\sum_{e}|h(e)|^2 < \infty.\tag6-19$$
In the last formula $e$ runs over all cells of $\tilde K$. Using the 
same arguments in the opposite way, we see that any cochain on
$\tilde K$ with values in $\C$ which satisfies $L^2$-condition comes from a
cochain with values in $\ell^2(\pi)$ considered as a $\C[\pi]$-module.
Comparing also the boundary maps shows that in this case the complex
$\om_{\Lambda}(C_\ast(\tilde K),\ell^2(\pi))$ can be identified with the
complex built on the $\C$-valued cochain with $L^2$ condition.

A slightly more general example is given by the module $V\otimes\ell^2(\pi)$,
where $V$ is a finite dimensional Hilbert representation of $\pi$. The
group $\pi$
acts diagonally on the tensor product above (which is understood over $\C$);
the corresponding von Neumann category is described in example 6 of 
section 2.6. This example may also be described geometrically similar to
the previous discussion as the cochain complex formed by chains
 $h$ on $\tilde K$ with values in $V$,
which satisfy the following $L^2$-condition 
$\sum_{e\subset \tilde K} |h(e)|^2 < \infty$.
The last sum it taken over all cells of $\tilde K$. 

\subheading{6.11} Here is a remark concerning building of $\C[\pi]$ modules
in the extended category. 

Generally, if we have two $\C[\pi]$-modules 
$\H_0$ and $\H_1$ in $\ca$ and a $\C[\pi]$-morphism $\alpha:\H_0\to \H_2$
then $\M=(\alpha:\H_0\to\H_1)$ represents a $\C[\pi]$-module in $\eca$.
Not all examples are of this form. 

For instance, consider the torsion object 
$\X_{\theta,\nu}=(\alpha:L^2(S^1)\to L^2(S^1))$ as in example 
3.12. Let the function $f\in L^\infty(S^1)$ be given by
$f(z)=\sqrt{1+\alpha(z)}$. Then the morphism
$m_f:\X_{\theta,\nu}\to \X_{\theta,\nu}$ 
(here $m_f$ denotes the operator of multiplication by $f$)
represented by the diagram 
$$
\CD
(\alpha: & \X_{\theta,\nu} & \to &  \X_{\theta,\nu})\\
   &        @V{m_f}VV              @VV{m_f}V\\
(\alpha: & \X_{\theta,\nu} & \to &  \X_{\theta,\nu})
\endCD
$$
defines an action of the group $\Z_2$ on $\X_{\theta,\nu}$.
Although it is only a $Z$ action on $L^2(S^1)$.

\subheading{6.12. Example: Extended cohomology of a mapping torus} For
illustrative purposes and because of many applications,
we will compute here the extended cohomology of a mapping torus.

Suppose that $Y$ is a finite polyhedron with a base point $y_0$ 
and let $\phi:Y\to Y$ be a homeomorphism with $\phi(y_0)=y_0$.
Consider the mapping torus $K\ =\ K_\phi \ =  Y\times [0,1]/\thicksim$, where 
$(y,1)\thicksim (\phi(y),0)$ for all $y\in Y$. We will consider $(y_0,0)$ as
the base point of $K$. 

We have the following exact sequence of the fundamental groups
$$0\to \pi(Y)\to \pi(K)\to \Z\to 0.\tag6-20$$
This sequence clearly splits. The natural splitting is given by sending 
the generator of $\Z$ to the class of the loop $(y_0,t)$, where $t\in[0,1]$;
the class of this loop in
$\pi(K)$ will be denoted $\tau$. We see that $\pi(K)$ is a semi-direct
product;  the following relation holds in $\pi(K)$:
$$\tau g \tau^{-1} = \phi_\ast(g)\tag6-21$$
for all $g\in\pi(Y)$, where $\phi_\ast:\pi(Y)\to\pi(Y)$ is the automorphism
induced by $\phi$. 

The universal covering $\tilde K$ of $K$ can be identified with the product
$\tilde Y\times \R$, where $\tilde Y$ is the universal covering of $Y$. 
To describe the action of the fundamental group $\pi(K)$
on $\tilde Y\times \R$
we have to fix a base point in $\tilde Y$ (above $y_0$) and then there is a
unique lift $\tilde \phi:\tilde Y\to \tilde Y$ of the homeomorphism $\phi$,
preserving the base points. 
Such lift $\tilde \phi$  {\it is equivariant} in the following sense:
$$\tilde \phi(g\cdot\tilde y) = \phi_\ast(g)\cdot \tilde\phi(\tilde y).\tag6-22$$
Any  element $g\in \pi(Y)$ acts on $\tilde K=\tilde Y\times \R$ 
by $g\cdot(\tilde y,t)=(g\cdot \tilde y,t)$ 
for all $\tilde y\in \tilde Y$, $t\in \R$, and the action of the element
$\tau$ is given by
$\tau\cdot(\tilde y,t)=(\tilde\phi(\tilde y),t-1)$.

The cell structure of $\tilde K$ is as follows. For any $n$-dimensional
cell $e\subset Y$, the universal covering $\tilde K$ has $n$-dimensional
cells of the form 
$\tau^i(ge)$, where $i\in\Z$, $g\in \pi(Y)$, and $(n+1)$-dimensional 
cells of the form
$\tau^i(ge\times (0,1))$, where $i\in\Z$, and $g\in \pi(Y)$. 

This allows to describe completely the chain complex $C_\ast(\tilde K)$ 
as a complex of left $\Lambda$-modules, where $\Lambda=\C[\pi(K)]$. Denote
$\Lambda_Y=\C[\pi(Y)]$, which will be understood as a subring of 
$\Lambda$ via (6-20).
 Then the morphism of complexes of left
$\Lambda$-modules
$$\tau^{-1}\otimes C_\ast(\tilde\phi): \Lambda\otimes_{\Lambda_Y}C_\ast(\tilde Y) \to 
\Lambda\otimes_{\Lambda_Y}C_\ast(\tilde Y),\tag6-23$$
which is given by the formula
$$\lambda\otimes c\mapsto \lambda\tau^{-1}\otimes C_\ast(\tilde \phi)(c),
\quad \lambda\in\Lambda,\quad c\in C_\ast(\tilde Y),\tag6-24$$
{\it is well defined, and is a chain morphism of complexes over $\Lambda$.}
Here $C_\ast(\phi): C_\ast(\tilde Y)\to C_\ast(\tilde Y)$ denotes 
the chain map
induced by $\tilde \phi$.
Comparing with the description of the cell structure of $\tilde K$ in the 
previous paragraph, one obtains that {\it the chain complex 
$C_\ast(\tilde K)$  coincides with the mapping cone of the following 
chain map}
$$\id - \tau^{-1}\otimes C_\ast(\tilde\phi): 
\Lambda\otimes_{\Lambda_Y}C_\ast(\tilde Y) \to 
\Lambda\otimes_{\Lambda_Y}C_\ast(\tilde Y).\tag6-25$$

Now, to compute the extended cohomology of $K$
(according to the basic construction of section 6.5) we have to apply the
functor $\om_\Lambda(\cdot,\M)$ to the mapping cone of
(6-25). Here $\M$ denotes
a $\Lambda$-module in the extended category $\eca$, which has to be chosen as 
the coefficient system. We obtain a representation of the complex 
$\om_\Lambda(C_\ast(\tilde K),\M)$ also as a mapping cone (with a slight
shift of the dimensions) of a morphism of
cochain complexes 
$$\om_{\Lambda_Y}(C_\ast(\tilde Y),\M|_Y) \to 
\om_{\Lambda_Y}(C_\ast(\tilde Y),\M|_Y)\tag6-26$$
in $\eca$. Here $\M|_Y$ denotes $\M$ but considered only as a left 
$\Lambda_Y$-module. 

Observe that defining $\M$ as a module over $\Lambda$
is equivalent to defining it as a module over $\Lambda_Y$ and specifying
an isomorphism of left $\Lambda_Y$-modules 
$$\tau: \M \to \phi^\ast\M,\tag6-27$$
which is given by multiplication by $\tau$;  
here $\phi^\ast\M$ denotes $\M$ with the $\Lambda_Y$-module structure 
on $\M$ given by $\rho_{\phi^\ast\M}=\rho_{\M}\circ\phi_\ast$, where
$\phi_\ast:\Lambda_Y\to\Lambda_Y$ is the automorphism induced by the
homeomorphism $\phi:Y\to Y$.
Thus, identifying the morphism (6-26), we obtain the following result:

\proclaim{6.12.1. Theorem} Let $K = K_\phi$ be the mapping torus 
of a homeomorphism $\phi:Y\to Y$ as above and let
$\M$ be a $\Lambda=\C[\pi]$-module in the extended category $\eca$,
where $\pi=\pi(K)$. Then there exists the following
exact sequence
$$
\aligned
\dots &\to \H^{i-1}(Y,\M)\to \H^i(K,\M)@>{i^\ast}>> \\
@>i^\ast>>
\H^i(Y,\M)
&@>{\id-\tau^{-1}\phi^\ast}>> \H^i(Y,\M)\to \H^{i+1}(K,\M)\to\dots
\endaligned\tag6-28
$$ 
of objects and morphisms of the extended category $\eca$.
In this exact sequence, 
$\phi^\ast: \H^i(Y,\M)$ $\to \H^i(Y,\phi^\ast\M)$ denotes the induced 
morphism (6-14), and
$\tau: \H^i(Y,\M)$ $\to \H^i(Y,\phi^\ast\M)$ denotes "the coefficient" 
morphism induced by isomorphism (6-27). \qed
\endproclaim

\subheading{6.12.2} Consider now a special case of the previous 
Theorem assuming that {\it the fundamental group $\pi(Y)$ acts 
trivially  on} $\M$. 

Then $\H^i(Y,\M)$ equals to $H^i(Y,\C)\otimes_\C \M$. We also observe that
the module $\phi^\ast\M$ is now identical with $\M$ and the induced map
$\phi^\ast: \H^i(Y,\M)\to\H^i(Y,\M)$ can be identified with the the 
tensor product of the induced map on the ordinary cohomology times the
identity on $\M$. 

We have given a left $\C[\tau, \tau^{-1}]$-module structure on some
$\M\in\ob(\eca)$. 
Also, we may view the ordinary homology and cohomology of $Y$ as left
modules over the ring $\C[\tau, \tau^{-1}]$ with $\tau$ acting by means of
the homomorphism $\phi_\ast$ induced by homeomorphism $\phi$. Having this
in mind, we will 
denote the kernel and cokernel of the following morphism of $\eca$
$$1\otimes\tau -\phi^\ast\otimes 1: H^i(Y)\otimes_\C \M \to 
H^i(Y)\otimes_\C\M\tag6-29$$
by
$$\Hom_{\C[\tau, \tau^{-1}]}(H_i(Y),\M)\quad\text{and}\quad
\Ext_{\C[\tau, \tau^{-1}]}(H_i(Y),\M)\tag6-30$$
correspondingly. Note that these notations can be justified.

Now, we may rewrite the exact sequence (6-28) in the 
following form
$$0\to \Ext_{\C[\tau, \tau^{-1}]}(H_{i-1}(Y),\M) \to \H^i(K,\M) \to
\Hom_{\C[\tau, \tau^{-1}]}(H_i(Y),\M) \to 0.\tag6-31$$

\subheading{6.12.3} In the situation of the previous subsection assume that
$\M\in\ob(\ca)$ (i.e. $\M$ is projective) 
and that the morphism $\tau:\M\to\M$ {\it has no discrete spectrum.} 

Then, as it is easy to see, $\Hom_{\C[\tau, \tau^{-1}]}(H_i(Y),\M)$ always
vanishes and we obtain the following result:

\proclaim{6.12.4. Corollary} Under the assumptions of Theorem 6.12.1,
suppose additionally that $\M\in\ob(\ca)$ is such that $\pi(Y)$ acts on it 
trivially and the isomorphism $\tau:\M\to\M$ has no discrete spectrum.
Then the extended cohomology of the mapping torus $K=K_\phi$ is given by
the following object of the extended category $\eca$
$$\H^i(K,\M) = (1\otimes\tau - \phi^\ast\times 1: H^{i-1}(Y)\otimes_\C\M \to
H^{i-1}(Y)\otimes_\C\M).\tag6-32$$
\endproclaim

\subheading{6.12.5. Example} Let $Z$ be a compact Hausdorff space and let
$\mu$ be a positive Radon measure on $Z$ such that the measure of any non-empty
open subset $U$ of $Z$ is positive. Consider the von Neumann category $\ca$
of representations of the algebra $L^\infty_\C(Z,\mu)$ of essentially 
bounded measurable functions on $Z$, which was discussed in section 2.6
(example 7). 

The simplest example of a $\C[\tau, \tau^{-1}]$-module in this category
is as follows. Let $\M$ be $L^2(Z,\mu)$ and let $\tau:\M\to \M$ be the 
operator of multiplication by a function $\tau:Z\to\C$, where
$\tau,\tau^{-1}\in L^\infty_\C(Z,\mu)$. 

Then the condition that $\tau$ has no discrete spectrum is equivalent to
the assumption that {\it all level sets $\{\xi\in\Z;\tau(\xi)=c\}$ have zero
measure.} 

Formula (6-32) gives an explicit representation of the extended cohomology.
The answer depends only on the structure of the Jordan decomposition of 
the induced morphism $\phi_\ast:H_{i-1}(Y) \to H_{i-1}(Y)$ on the ordinary 
homology with complex coefficients.

In particular, using Theorem 4.10 we obtain that {\it if the function
$\tau$ is continuous and all the level sets have measure zero,
then the $i$-dimensional 
extended cohomology $\H^i(K,\M)$ of the mapping torus $K$ is nonzero if 
and only if the spectrum of the induced map 
$\phi_\ast:H_{i-1}(Y)\to H_{i-1}(Y)$ 
has nonempty intersection with the set of values of function $\tau$,
$\tau(Z)\subset \C$.}

As a conclusion we obtain:
 
{\it if the induced map 
$\phi_\ast:H_{i-1}(Y)\to
H_{i-1}(Y)$ is semi-simple, then studying extended cohomology of the
mapping torus, one may recover all information about 
the usual homology of $Y$ with complex coefficients and the action
of the homeomorphism $\phi_\ast$ on it}.  

This should be compared with the fact that the reduced $L^2$ cohomology
of the mapping torus vanishes in many cases, cf. \cite{L1}.

Note also that using Theorem 6.12.1 
inductively one may hopefully 
calculate the extended cohomology of abelian groups and 
of nilpotent groups.

\heading {\bf \S 7. De Rham theorem for extended cohomology}\endheading

The purpose of this section is to define De Rham version of the extended
$L^2$ cohomology and to prove a De Rham type theorem, establishing
isomorphism between combinatorially and analytically defined objects.

It should be noticed that some important consequences of this theorem
are already well known: they are theorem of J.Dodziuk \cite{D}
(which gives equality of the von Neumann dimensions) and theorem of 
A. Efremov \cite{E}, \cite{E1} (which establishes equivalence of the
spectral density functions and equality of the 
Novikov-Shubin invariants). Since we know that the isomorphism type of a 
torsion object is not determined by the spectral density function, there is
still a need to compare fully the combinatorially and analytically defined
objects of the extended abelian category. 

To my mind, the main part of \cite{E}, \cite{E1} 
are written without control on the category, where the arguments 
are performed. This creates difficulties while trying to adopt the
result of \cite{E}, \cite{E1} in the framework of the present approach. 
As an example, let me 
mention that admitting densely defined operators, as morphisms (as in
\cite{E}, \cite{E1}) immediately
makes all injective with dense image operators invertible and so 
completely eliminates the torsion phenomenon. I do not exclude possibility 
that the arguments of \cite{E}, \cite{E1} could be rewritten and
presented in a different
form to avoid the categorical difficulties. Instead of doing this, I suggest 
in this section arguments (using sheaves and spectral sequences) which are
much easier and more transparent, than the 
original arguments of \cite{D} and \cite{E}, \cite{E1}.

I should mention that a new preprint of
M. Shubin \cite{S} contains a different version of De Rham type theorem
for extended $L^2$ cohomology. 

A few words about notations.
In this section $\A$ will denote a fixed algebra with involution and we 
will denote by $\ca$ the 
von Neumann category of all Hilbertian representations of $\A$, cf. \S 5. 
We will need also a category containing $\ca$, which we will denote
by $\fa$. Its objects are all Fr\'echet topological vector spaces supplied 
with an action of the algebra $\A$ (no condition involving the involution)
and the morphisms are continuous $\A$-linear 
maps. 

\subheading{\bf 7.1. The De Rham complex} First, we are going to discuss the
notion of {\it flat bundle} with fiber a Hilbertian
module. This notion is quite similar to the 
standard notion of flat finite dimensional vector bundle, or to the notion
of a flat Hilbert bundle over a von Neumann algebra, which was considered
in \cite{BFKM}, for example; cf. also \cite{CFM}. Roughly, flat bundle is
a bundle with discrete structure group.
                                                                    
Let $\M$ denote a fixed Hilbertian representation of $\A$. 
This module $\M$ will serve as the typical fiber of our flat bundles. 
$GL(\M)$ will denote the group of all $\ca$-automorphisms of $\M$, i.e. all
linear homeomorphisms $\M\to \M$ commuting with the $\A$-action. 

Let $X$ be a topological space. {\it A flat Hilbertian bundle over $X$ 
with fiber} 
$\M$ can be defined by a 1-cocycle on an open covering $\U$ of $X$ 
with values in the group $GL(\M)$, i.e.
by a function which associates an element $g_{UV}\in GL(\M)$ for any ordered
pair of open sets $U,V\in \U$ having nonempty intersection. It satisfies the 
following conditions: (a) $g_{UU}=1$ and 
(b) $g_{UV}\cdot g_{VW}=g_{UW}$ if the
intersection of the sets $U$, $V$ and $W$ is nonempty. As usual, given a
1-cocycle as above, one can construct the space of the bundle
$\E$ identifying in the disjoint union of the spaces $U\times \M$ the points
$(x,m)\in U\times \M$ with $(x, g_{UV}\cdot x)\in V\times \M$ for $x$ 
belonging to the intersection of $U$ and $V$. The space $\E$ admits a 
natural projection onto $X$. The flat Hilbertian bundle as above will be 
denoted by one letter $\E$. There is a well-known equivalence relation 
between the 1-cocycles.

If $X^\prime\subset X$ is a subset, there is naturally
defined the restriction of a flat Hilbertian bundle $\E$ over $X$, onto $X'$,
which is denoted\ $\E|_{X^\prime}$ and is a flat Hilbertian
bundle over $X^\prime$.

Suppose now that $X$ is a finite-dimensional smooth manifold 
(not necessarily compact).
Denote by $C^\infty_\M(X)$ the Fr\'echet space of the smooth functions on $X$ 
with values in $\M$. The family of Fr\'echet spaces $C^\infty_\M(U)$, where 
$U$ runs over all open subsets of $X$ together with the natural restriction 
maps, form a sheaf with values in the category $\fa$, 
which will be denoted by $C^\infty_\M$.

We want to study {\it smooth differential forms on $X$ with
values in the flat Hilbertian bundle} $\E$ over $X$, cf. above. 
If $\E$ is represented by a 1-cocycle with respect to a finite open covering 
$\U=\{U\}$ of $X$, then for each open set $U\in\U$ define $A^i(U,\E)$ as 
$$A^i(U,\E)\ =\ \Lambda^i(U)\otimes C^\infty_\M(U),\tag7-1$$
where $\Lambda^i(U)$ denotes the space of smooth $i$-forms on $U$ with 
complex values, and the tensor product is taken over the ring $C^\infty(U)$. 
{\it A differential $i$-form on $X$ with values in $\E$} is defined as a 
collection
$$\omega\ =\  =\ \{\omega_U\}_{U\in\U},\tag7-2$$
where $\omega_U\in A^i(U,\E)$ such that for any pair $U, V\in\U$ with
nonzero intersection we have
$$\omega_U|_{U\cap V}\ =\ g_{UV}\cdot \omega_V|_{U\cap V}.\tag7-3$$
The space of all smooth $i$ forms on $X$ with values in $\E$ will be denoted
$A^i(X,\E)$. This space has a natural Fr\'echet topology. 
Also, the initial
algebra $\A$ acts on $A^i(X,\E)$.

There is {\it the operator of covariant derivative}
$$\nabla : A^i(X,\E)  \to A^{i+1}(X,\E),\tag7-4$$ 
defined as follows. If $\omega\in A^i(X,\E)$ is given by 
$\omega\ =\  =\  \{\omega_U\}_{U\in\U}$, then 
$\nabla(\omega)\ =\  \{d\omega_U\}_{U\in\U}$, where $d$ denotes the usual 
exterior derivative. The operator $\nabla$ is clearly well defined, 
continuous, commutes with the $\A$-action and
satisfies $\nabla^2\ =\ 0$.

We obtain the cochain complex in $\fa$, the {\it De Rham complex with values 
in $\E$}:
$$0\to A^0(X,\E)@>\nabla>> A^1(X,\E)@>\nabla>>\dots 
A^i(X,\E)@>\nabla>>
A^{i+1}(X,\E)@>\nabla>> \dots A^n(X,\E)\to 0,\tag7-5$$
where $n$ is the dimension of the manifold $X$. 
We will denote the De Rham complex by $(A^\ast(X,\E),\nabla)$ for short.

\subheading{7.2. Extended $L^2$ \v Cech cohomology} Suppose that $\E$ is a 
flat Hilbertian bundle over a closed manifold $X$. Let $\U$ be a {\it good} 
finite open covering
of $X$ (cf. \cite{BT}, page 42) and suppose that the flat bundle $\E$ is
trivial over every open set $U\in \U$. We will view now the flat bundle
$\E$ as a sheaf of its flat sections; it is defined as follows
$$\Gamma_\E(W)\ =\ \ker [\nabla:A^0(W,\E)\to A^1(W,\E)]
\quad\text{for any
open subset}\quad W\subset X.\tag7-6$$
Note that $\Gamma_\E(W)$ is isomorphic to $\M$ if $\E|_W$ is trivial.
We obtain in this way the space of \v Cech cochains $C^p(\U;\E)$ (defined in 
the standard way) and the
\v Cech cochain complex 
$$\dots @>\delta>>C^p(\U;\E)@>\delta>>C^{p+1}(\U,\E)@>\delta>>\dots\tag7-7$$
which we will denote by $(C^\ast(\U;\E),\delta)$. 
Observe that 
the spaces of cochains $C^p(\U;\E)$ are naturally defined as 
{\it Hilbertian modules and the \v Cech boundary homomorphisms are
continuous linear maps}. Thus, the \v Cech complex can be viewed as a
complex in category $\ca$. Moreover, if the fiber $\M$ is finite, i.e.
$\M\in\ob(\caf)$, then the whole \v Cech complex $C^p(\U;\E)$ belongs
to $\caf$. 

In particular, it follows that we may define {\it the 
\v Cech version of the extended $L^2$ cohomology} (which we will denote by
$\check \H^\ast(\U,\E)$) as the cohomology (understood in $\E(\ca)$ or in 
$\E(\caf)$) of the above \v Cech complex. 
It follows from Theorem 7.3 below that the \v Cech extended
cohomology $\check \H^\ast(\U,\E)$ does not depend on the (finite good open)  
cover $\U$.
                  
Note that the construction of the extended $L^2$ cohomology used in \S 6
(based upon cell decompositions) 
is a special case of this \v Cech construction. Namely, choose a sufficiently
fine $C^1$ tringulation of the manifold $X$ and consider the covering 
$\U$ of $X$ formed by the stars of the vertices. Then it is easy to see that
the \v Cech cohomology $\check \H^\ast(\U,\E)$ of this covering
will coincide with $\H^\ast(X,\M)$, defined as in \S 6, using the monodromy
representation of the flat bundle $\E$.

\proclaim{7.3. Theorem. (\v Cech - De Rham Theorem)} Suppose that $X$ is a 
closed smooth manifold
and $\E$ is a flat Hilbertian bundle over $X$ with fiber a Hilbertian module
$\M$ over $\A$; let $((A^\ast(X,\E),\nabla)$ be the corresponding
De Rham complex. 
Suppose that $\U$ is a finite good (cf. \cite{BT}) open 
cover of $X$ and let $C^\ast(\U,\E)$ be the corresponding \v Cech complex, 
viewed as a complex in category $\ca$ of Hilbertian modules. 
Then there is a chain homotopy equivalence 
$$(C^\ast(\U,\E),\delta)\to (A^\ast(X,\E),\nabla)\tag7-8$$ 
in $\fa$. 
\endproclaim

In particular, we obtain, that the \v Cech extended $L^2$ cohomology
$\H^\ast(\U,\E)$ does not depend on the choice of a covering $\U$. We will
mention some other corollaries later, cf. 7.10, 7.11.

The proof of Theorem 7.3 (which is given in 7.7)
will be quite similar to the proof of the 
standard De Rhan theorem given in \cite{BT} for example. The main 
difference will be in the fact that we are not allowed to use spectral 
sequences in category $\fa$ (since it is only additive and not abelian). 
Instead, we will use Lemma 7.6. 

The main step is the Poincar\'e lemma.

\proclaim{7.4. Lemma. (Poincar\'e lemma)} Let $X$ be diffeomorphic to $\R^n$.
Let $M$ be a finitely generated Hilbertian module and let $\E$ denote the
trivial flat Hilbertian bundle with fiber $M$ over $X$. Then there is a
homotopy equivalence 
$$F: (\M, 0) \ \to\ (A^\ast(X,\E),\nabla)\tag7-9$$
between cochain complexes in category $\fa$. Here $(\M, 0)$ denotes the 
cochain complex, which has only one nonzero chain space $\M$ in dimension
zero. The map $F:\M\to A^0(X)$ is given by assigning to a point 
$m\in\M$ the constant section with value $m$.
\endproclaim

\demo{Proof} The proof of the Poincar\'e Lemma repeats the arguments 
of \cite{BT},
pages 34, 35 with an additional attention to the fact that the constructed 
in \cite{BT} standard chain homotopies are actually morphisms of category 
$\fa$. \qed
\enddemo

\subheading{\bf 7.5. Double complexes in $\fa$} For the proof of Theorem 7.3
we will need some properties
of double complexes in category $\fa$. They have in fact a very general
nature and are valued in arbitrary additive categories. In abelian categories
they trivially follow from existence of two spectral sequences of a double
complex.

Consider a {\it double complex} $K\ =\ \oplus_{p,q} K^{pq}$ in 
category $\fa$ having two differentials $D^\prime$, of bidegree $(1,0)$, and 
$D^{\prime\prime}$,
of bidegree $(0,1)$, which satisfy ${D^\prime}^2=0$, ${D^{\prime\prime}}^2=0$
and $D^\prime D^{\prime\prime}+ D^{\prime\prime}D^{\prime}=0$.
The {\it total differential} $D\ =\ D^\prime\ + \ D^{\prime\prime}$ 
then satisfies $D^2=0$. We will always assume that all our double complexes
are {\it positive and bounded} which means that there exists some integer
$N$ such that $K^{pq}$ vanishes if one of $p$ and $q$ is negative or larger
than $N$. Because of this condition we may consider {\it the total complex}
$\Tot(K)\ =\ K^n$, where $K^n=\oplus _{p+q=n}K^{pq}$, considered with 
the total differential $D$. 

If $K$ and $L$ are two such double complexes, {\it a morphism} $\phi:K\to L$
is a collection of morphisms $K^{pq}\to L^{pq}$ in $\fa$, which commute
with the differentials $D^\prime$ and $D^{\prime\prime}$. Morphism of 
double complexes clearly determines a chain map of the total complexes.
We are interested in conditions on $\phi$ which guarantee that this chain
map of the total complexes is a chain homotopy equivalence. 

\proclaim{7.6. Lemma} Let $K$ and $L$ be positive and bounded
double complexes in $\fa$. Suppose that $\phi:K\to L$ is a morphism between
these double complexes satisfying the following condition: for each integer
$p$ the induced chain map of the {\it vertical} complexes
$(K^{p,\bullet},D^{\prime\prime})\to (L^{p,\bullet},D^{\prime\prime})$
is a homotopy equivalence. Here $(K^{p,\bullet}\ =\ \oplus _{q}K^{pq}$ and
similarly for $L$). Then the induced by $\phi$ morphism of the total
chain complexes $(\Tot(K),D)\to (\Tot(L),D)$ is a chain homotopy equivalence.

The same conclusion is true if we assume that for every $q$ the induced 
by $\phi$ chain map of the {\it horizontal} complexes
$(K^{\bullet q},D^{\prime})\to (L^{\bullet q},D^{\prime})$ 
is a homotopy equivalence.
\endproclaim
\demo{Proof} Suppose first that our double complexes have only two nontrivial 
columns, i.e. $K^{pq}=0$ and $L^{pq}=0$ for all $p\ne i,\ i+1$. Then 
the double complex $K$ consists of two cochain complexes $K^{i,\bullet}$ 
and $K^{i+1,\bullet}$ (two columns) and the horizontal differential 
$D': K^{i,\bullet} \to K^{i+1,\bullet}$ can be viewed as a chain map
between them. Then the total complex
coincides with the cone of this chain map; similarly for $L$. In this 
form Lemma 7.6 is well known: it is one of the basic properties
of triangulated categories; cf. for example \cite{KS}, Corollary 1.5.5.

The general case can be reduced to the case of two columns using induction. 
Namely, suppose $K^{i,\bullet}$ and $L^{i,\bullet}$ are the last nontrivial
columns, i.e. $K^{pq}=0$ and $L^{pq}=0$ for $p>i$. Let $K'$ be the total
complex of the truncated double complex $\oplus_{p<i}K^{pq}$; similarly
we will construct $L'$. Then by induction we know that the induced by $\phi$
chain map $K'\to L'$ is a chain homotopy equivalence. The horizontal
differential $D'$ gives us the maps $D':K'\to K^{i,\bullet}$ and  
$D':L'\to L^{i,\bullet}$ and the cones of these maps coincide with the
total complexes $\Tot(K)$ and $\Tot(L)$ correspondingly. Thus we may again
apply the Lemma in the case of two columns and the result follows. \qed
\enddemo
\subheading{\bf 7.7. Proof of theorem 7.3} For any open set $U\subset X$ we 
may consider 
the space of smooth $q$-forms on $U$ with values in $\E$, which we denote now
$A^q_\E(U)$; this defines a sheaf $A^q_\E$ over $X$ with values in category 
$\fa$.
The operator of covariant derivative defines a homomorphism of sheaves
$\nabla: A^q_\E\to A^{q+1}_\E$.   

Fix a good finite open cover $\U$ of $X$.

Consider the following double complex $K=\oplus_{pq}K^{pq}$, where
$K^{pq}\ =\ C^p(\U,A^q_\E)$ is the space of \v Cech cochains of cover $\U$
with values in sheaf $\A^q_\E$; this construction was called in \cite{BT}
the \v Cech - De Rham complex. We have two differentials in this double 
complex: $D^\prime=\delta$ is the \v Cech differential, and 
$D^{\prime\prime}=(-1)^p\nabla$ on $K^{pq}$ is the covariant derivative. 
These two differentials satisfy the necessary commutation relations. Note that
the constructed double complex is positive and bounded. 

We will construct now two more double complexes. The first one, which we will
denote by $A$, is obtained from the De Rham complex. Namely, define
$A^{pq}$ to be $0$ if $p\ne 0$ and $A^{0q}=A^q_\E(X)$. The horizontal differential
$D^\prime$ will be zero and the vertical $D^{\prime\prime}$ coincides with
$\nabla$.

The third double complex will be denoted by $\Check C$, it is essentially the
\v Cech complex. It is defined by $\Check C^{pq}=0$ for $q\ne 0$ and 
$\Check C^{p0}=C^p(\U,\E)$. Note that $\E$ is considered here as the sheaf
of its flat sections. The horizontal differential of $\Check C$ is the 
\v Cech boundary $\delta$ and the vertical differential is zero.

There are two natural morphisms of the constructed double complexes
$$\phi: A\ \to \ K\quad\text{and}\quad \psi: \check C\ \to\ K$$
given as follows: $\phi^{0q}:A^{0q}\to C^0(\U,A^q_\E)$ acts by 
$\phi^{0q}(\omega)=ga)=\{\omega|_U\}_{U\in\U}$ and the morphisms on all other
places are zero. The second morphism $\psi$ is the natural inclusion
$\psi^{p0}:C^{p0}\ =\ C^p(\U;\E)\to C^p(\U;A^0_\E)$ and zero on all other
places.

Our theorem 7.3 will follow from Lemma 7.6 (applied twice)
if we will establish two the following statements:
 
(1) for each integer $q$ the induced horizontal chain map 
$\phi^{\bullet q}: A^{\bullet q}\to K^{\bullet q}$ is a chain homotopy
equivalence in $\fa$;

(2) for each integer $p$ the induced  vertical chain map
$\psi^{p\bullet}:\Check C^{p\bullet}\to K^{p\bullet}$ is a chain homotopy
equivalence in $\fa$.

However, (2) is guaranteed by the Poincar\'e lemma 7.4 and (1) is given by 
the chain homotopy constructed in \cite{BT}, page 94 using a partition of 
unity. This completes the proof. \qed

\subheading{7.8. De Rham version of the extended cohomology} Unfortunately,
one cannot use the De Rham complex (7-5) in order to define the De Rham
version of the extended cohomology. The reason is that (7-5) is a complex
in category $\fa$ of Fr\'echet representations of $\A$, 
and the abelian extension constructions of \S 1 and \S 5 do not apply
to it. In order to overcome this difficulty, we consider here another
complex (the Sobolev - De Rham complex) 
which is homotopy equivalent to (7-5) and belongs to category $\ca$ of
all Hilbertian representations of $\A$; note that $\ca$ admits an extended 
abelian category $\eca$, which was constructed in \S 5.

If $X$ is a closed smooth manifold and $\E$ is a flat Hilbertian bundle 
over $X$, we will denote by $A^p_{(m)}(X,\E)$ the $m$-th Sobolev space
of $p$-forms on $X$ with values in $\E$ (here $m$ is a non-negative integer). 
This space is defined as the completion of the
space of smooth forms $A^p(X,\E)$ (which appeared earlier) with respect
to the $m$-th Sobolev's norm $||\ ||_{(m)}$. The latter is defined in the
usual way by using partitions of unity.

The covariant derivative $\nabla$ defines a continuous operator
$\nabla: A^p_{(m)}(X,\E) \to A^{p+1}_{(m-1)}(X,\E)$. Thus we obtain {\it the
$m$-th Sobolev - De Rham complex}
$$0\to A^0_{(m)}(X,\E)\dots \to A^{p-1}_{(m+1-p)}(X,\E)@>\nabla>> 
A^p_{(m-p)}(X,\E) @>\nabla>> A^{p+1}_{(m-1-p)}(X,\E)\to \dots\tag7-10$$
Here $m\ge n=\dim X$. 

The complex (7-10) was studied earlier in \cite{D} and in \cite{LL} in
special cases. 

Note that (7-10) is a
chain complex in category $\ca$, i.e. it consists of Hilbertian spaces
(cf. \cite{P}) and continuous linear maps commuting with $\A$. Thus, 
using the results of \S 5, we may
define {\it the De Rham version of the extended $L^2$ cohomology} 
$$\H^p_{(m-p)}(K,\E) = (\nabla: A^{p-1}_{(m+1-p)}(X,\E)\to Z^p_{(m-p)}),
\tag7-11$$
where $Z^p_{(m-p)} = \ker[\nabla: A^p_{(m-p)}(X,\E)\to 
A^{p+1}_{(m-p-1)}(X,\E)],$
i.e. the cohomology of (7-10) in the extended category $\eca$.

The following theorem implies, in particular, that the extended cohomology 
of the Sobolev - De Rham complex (7-10) does not depend on $m$.

\proclaim{7.9. Theorem. (Sobolev - De Rham Theorem)} Let $X$ be a closed
smooth manifold with a flat Hilbertian bundle $\E\to X$. Then the $m$-th
Sobolev - De Rham
complex (7-10) is homotopy equivalent in category $\ca$ to the \v Cech complex 
$(C^p(\U;\E),\delta)$ constructed using any finite good open cover $\U$
of $X$.
\endproclaim
\demo{Proof} We will show that the natural inclusions 
$i: A^p(X,\E) \to A^p_{(m-p)}(X,\E)$ determine a chain homotopy equivalence
$i$ in $\fa$ between the De Rham complex (7-5) and the $m$-th 
Sobolev - De Rham complex (7-10); this together with Theorem 7.3 will complete
the proof of Theorem 7.9.

Fix a Riemannian metric on $X$ and a Hermitian metric on $\E$. Let
$\Delta_p = \nabla^\ast\nabla+\nabla\nabla^\ast$ denote the Laplassian 
acting on the spaces of $p$-forms. It determines a
morphism $\Delta_p: A^p_{(m-p)}(X,\E) \to A^p_{(m-p-2)}(X,\E)$ in $\ca$
for any 
$m$. Consider the heat operator
$e^{-\Delta_p}$. It is an infinitely smoothing operator, it maps 
any Sobolev space $A^p_{(m-p)}(X,\E)$ into the space of smooth forms
$A^p(X,\E)$, cf. for example, \cite{BFKM}.

We denote by $j: A^p_{(m-p)}(X,\E) \to A^p(X,\E)$ the map given by
$j(\omega) = e^{-\Delta_p}(\omega)$. It commutes with the covariant
derivative $\nabla$ 
and determines a chain map in category $\fa$ between the 
Sobolev and smooth De Rham complexes.  

We want to show that $i$ and $j$ are mutually inverse homotopy equivalences. 
To do so it is enough to construct an operator $H$ which maps continuously 
any Sobolev space $A^p_{(m')}(X,\E)$ into $A^{p-1}_{(m'+1)}(X,\E)$ and
satisfies the homotopy relation
$$j(\omega) - \omega = \nabla H + H\nabla,\quad\text{for any}\quad
\omega\in A^p_{(m')}(X,\E).\tag7-12$$

We define $H$ as follows
$$ H = H'\circ \nabla^\ast,\qquad\text{where}\quad
H' = \frac{1}{2\pi i}\int_C\frac{\phi(\lambda)}
{\lambda - \nabla^\ast\nabla}d\lambda,\tag7-13$$
and
$$\phi(\lambda) = \frac{e^{-\lambda}-1}{\lambda}.\tag7-14$$
Here we view $\nabla^\ast\nabla$ as a self adjoint operator acting from the
closure of the subspace of coclosed forms 
$\overline{\im(\nabla^\ast)}\subset A^p_{(m-p)}(X,\E)$ into itself, 
and being zero on the space orthogonal to $\overline{\im(\nabla^\ast)}$.  
It's spectrum is real and non-negative and so the
operator calculus applies. The integration curve $C$ is $x+1=|y|$ 
on the plane $(x,y)$ with the orientation chosen such that $y$ decreases. 
The function $\phi(\lambda)$ decays as $x^{-1}$
along the curve $C$ and the resolvent $(\lambda - \nabla^\ast\nabla)^{-1}$
also decays as $x^{-1}$; this shows convergence of $H$ in
(7-13) as an operator from $A^p_{(m')}(X,\E)$ to $A^p_{(m'+ 2)}(X,\E)$.

The homotopy relation (7-12) can be easily verified. \qed
\enddemo

We refer to the recent preprint \cite{S} of M. Shubin which contain a different version
of De Rham type theorem for extended cohomology.

\proclaim{7.10. Corollary} {\bf (J.Dodziuk \cite{D}, A.Efremov \cite{E})} 
Let $\ca^{fin}$ denote a full
finite von Neumann subcategory of $\ca$, cf. \S 5. Suppose that the fiber 
$\M$ of the flat bundle
$\E$ is finite, i.e $\M\in\ob(\ca^{fin})$. 
Let $\tr$ be
a trace on the category $\ca^{fin}$, cf. \S\S 2, 5. Then 
the von Neumann dimension with respect to the trace $\tr$ of the 
extended De Rham cohomology $\H^i_{(m-i)}(K,\E)$ coincides with the von 
Neumann dimension of the extended \v Cech cohomology $\check \H^i(\U,\E)$.
Also, the spectral density functions of the torsion parts of 
$\H^i_{(m-i)}(K,\E)$ and $\check \H^i(\U,\E)$ are dilatationally equivalent 
and their Novikov-Shubin invariants coincide.
\endproclaim

\proclaim{7.11. Corollary} {\bf (Finiteness.)} As in 7.10, suppose that 
given a full finite von Neumann subcategory 
$\ca^{fin}$ of $\ca$ and that $\M\in\ob(\ca^{fin})$. 
Then the extended De Rham cohomology $\H^i_{(m-i)}(K,\E)$ is also finite 
(i.e. it isomorphic in $\eca$ to an object of $\E(\ca^{fin}))$.
This can also be stated by saying that the Sobolev - De Rham complex 
$(A^\ast_{(m-\ast)}(K,\E),\nabla)$ is Fredholm.
\endproclaim

This follows automatically from Theorem 7.9 since the \v Cech
extended cohomology is clearly finite if $\M$ is finite.

\enddocument